\documentstyle[12pt]{article}
\newcommand{\noin}{\noindent}

\begin{document}

\title {On isometric reflexions in Banach spaces}

\bigskip

\author {A. Skorik and M. Zaidenberg}

\bigskip

\date{}
\maketitle

\hspace{3.5in} Dedicated to Professor S. G. Krein

\bigskip

\begin{abstract}
\noin We obtain the following characterization of Hilbert spaces. Let $E$ be a
Banach space whose unit sphere $S$ has a hyperplane of symmetry. Then $E$ is a
Hilbert space iff any of the following two conditions is fulfilled:

\noin a) the isometry group ${\rm Iso}\, E$ of $E$ has a dense orbit in S;

\noin b) the identity component $G_0$ of the group ${\rm Iso}\, E$ endowed with
the strong operator topology acts topologically irreducible on $E$.

Some related results on infinite dimentional Coxeter groups generated by
isometric reflexions are given which allow to analyse the structure of
isometry groups containing sufficiently many reflexions.

\end{abstract}

\bigskip

\section*{INTRODUCTION}

Let $E$ be a real Banach space, $S=S(E)$ the unit sphere in $E$, ${\rm Iso}\,
E$ the isometry group of $E$ endowed with the strong operator topology, and
$G_0 = G_0(E)$ the identity component of ${\rm Iso}\, E$. {\it A reflexion} in
$E$ is an operator of the form $s_{e, e^*} = 1_E - 2e^* \otimes e $, where
$e\in E, e^* \in E^*$ and $e^*(e)=1$. If $s=s_{e, e^*} \in {\rm Iso}\, E$, then
one may assume also that $ ||e||_E = ||e^*||_{E^*} =1$; in this case we will
call $e$ {\it the reflexion vector} and $e^*$ {\it the reflexion functional};
regarding as sphere points, $e$ and $-e$ are called {\it reflexion points}. The
unit sphere $S$ is symmetric with respect to {\it the mirror hyperplane} ${\rm
Ker}\,e^*$ of $s$. It turns out that this imposes strong restrictions on the
isometry group ${\rm Iso}\, E$.

We say that a proper subspace $H \subset E$ is {\it biorthogonally complemented
in} $E$ if there exists a bicontractive projection $p$ of $E$ onto $H$, i.e.
such that $ ||p||_E = ||1_E - p||_E = 1$.\bigskip

\noin {\bf Theorem 1.} {\sl  Let $s_{e, e^*}$ be an isometric reflection in
$E$. Let $H = {\overline{\rm span}}(G_0\, e)$ be the minimal closed subspace of
$E$ containing the orbit $G_0\, e$. Then

\smallskip

\noin a) $H$ is a Hilbert space and $H$ is biorthogonally complemented in $E$
or $H = E$;

\smallskip

\noin b) furthermore, there exists a projection $p$ of E onto $H$ such that

\smallskip

\noin i) $1_E - 2p \in {\rm Iso}\, E$,

\noin ii) $(1_E - p) + {\bar u} p \in {\rm Iso}\, E$ for any ${\bar u} \in {\rm
O}(H) = {\rm Iso} H$, and

\noin iii) any $g \in {\rm Iso}\, E$ such that $g \vert H = {\bar u} \in {\rm
O}(H)$ has the form $g= {\bar v} (1_E - p) + {\bar u} p$, where ${\bar v} = g
\vert {\rm Ker}\, p \in {\rm Iso}\, {\rm Ker}\, p$;

\smallskip

\noin c) the orbit $G_0\, e$ coincides with the unit sphere $S(H)$ of $H$.
}\bigskip

This subject is related to the following Banach - Mazur rotation problem ([3],
p.242):

\smallskip

{\sl Let $E$ be a separable Banach space such that the group ${\rm Iso}\, E$
acts transitively on the unit sphere $S$. Is it true that $E$ is a Hilbert
space?}

\smallskip

Recall (see [23, Ch.IX], \S 6) that the group ${\rm Iso}\, L_p$, where $L_p =
L_p [0;1]$ and $1 \le p \ne 2 < \infty$, has exactly two orbits on the unit
sphere $S_p = S(L_p)$. One of them consists of the functions in $S_p$ with the
zero set of a positive measure, and the other one contains the rest. Thus, both
orbits are dense in $S_p$. One says that the group ${\rm Iso}\, E$ acts {\sl
almost transitively} on $S$ if it has a dense orbit in $S$. This is the case in
the above examples and also in the anisotropic spaces $L_{pq}$. In a
non-separable $L_p-$space the second of the above two orbits is empty, and thus
 it  is a non-Hilbert Banach space with the isometry group acting transitively
on the unit sphere. This shows that the assumption of separability in the
Banach - Mazur problem is essential.

Observe that ${\rm Iso}\, E$ is a Banach-Lie group. If this group is transitive
on the unit sphere $S$, then $S$ is a homogeneuos space of ${\rm Iso}\, E$. If
in addition $S$ has a hyperplane of symmetry $L$, it should be a symmetric
space. Indeed, $L$ is a mirror hyperplane of an isometric reflection. The unit
sphere $S$ having a reflexion point, by transitivity each point $x \in S$
should be a reflexion point of an isometric reflexion $s = s_{x, x^*}$.
Furthermore, $x$ is an isolated fixed point of the involution $-s\,|\,S$ which
acts as $-1$ at the supporting hyperplane $x^* = 1$ to $S$ at $x$ (we are
grateful to J. Arazy for this remark). From Theorem 2 below it follows that $S$
being a symmetric space of the group ${\rm Iso}\, E$,  $E$ should be a Hilbert
space. In fact, in Theorem 2 more strong criteria for $E$ to be a Hilbert space
are done. They hold without the separability assumption.

\bigskip

\noin { \bf Theorem 2.} {\sl Let the group ${\rm Iso}\, E$ contains a
reflection $s_{e, e^*}$ along the vector $e \in S$. Then $E$ is a Hilbert space
iff either of the following two conditions is fulfilled:

\noin a) ${\rm Iso}\, E$ acts almost transitively on $S$;

\noin b) $e$ is a cyclic vector of the strong identity component $G_0$ of ${\rm
Iso}\, E$ (i.e. $E = {\overline {\rm span}} (G_0\, e)$).}  \bigskip

The second statement is a corollary of Theorem 1; the first one, being much
simpler, is proven along the same lines.

By a theorem of Godement [9] any isometric operator in a Banach space has a
non-trivial invariant subspace (see also [28] for a more general fact). From
Theorem 1 one obtains the following \bigskip

\noin { \bf Corollary.} {\sl Let $E$ be a non-Hilbert Banach space. If there is
an isometric reflexion in $E$, then all operators in $G_0(E)$ have a common
non-trivial invariant Hilbert subspace $H$, biorthogonally complemented in $E$.
Moreover, if $G_0(E)$ is a non-trivial group, then ${\rm dim}\, H > 1$. In
particular, in this case there is an orthogonally complemented euclidean plane
in $E$.}  \bigskip

Note that by a theorem of Yu. Lyubich [20] if a finite dimensional Banach space
has an infinite isometry group, i.e. if the group $G_0(E)$ is non-trivial, then
$E$ has a euclidean plane $L$ with a contractive projection $p : E \to L$ (in
this case $L$ is called {\it orthogonally complemented} in $E$) (see also [16],
[21]). From the other hand, there are Banach spaces of infinite dimension with
big isometry groups, but without any orthogonally complemented euclidean
subspace of dimension greater than $1$. Indeed, $L_p = L_p [0;1]$, where $1<p
\ne 2 < \infty$, contains no such a subspace, while the group $G_0$ is
non-trivial. Furthermore,  there is no bicontractive projection of $L_p\,\, (p
\ne 2)$ onto a hyperplane [13, 14]; in particular, there is no isometric
reflexion. The same is true in general for rearrangement-invariant (r.i.) ideal
Banach latticies, or symmetric spaces, of (classes of) measurable functions
different from $L_2$ [14, Theorem 4.4]. Recall [17, 19] that a r.i. (or
symmetric) space $E$ on the interval [0;1] satisfies the following axioms:

1) $1 \in E$ and $ ||1||_E =1$.

2) For any measure preserving transformation $\alpha$ of the interval [0;1] the
{\it shift operator} $T_{\alpha} : x(t) \rightarrow x(\alpha(t))$ acts
isometrically in $E$.

3) If $x(t) \in E$ and $ |y(t)| \le  |x(t)| $ a.e., then $y(t) \in E$ and $
||y(t)||_E \le  ||x(t)||_E$.

If $E$ is a r.i. space different from $L_2$, then every $g \in {\rm Iso}\, E$
has {\it a weighted shift representation} $g: x(t) \rightarrow h(t)x(\phi(t))$,
where $h = g(1) \in E$ and $\phi$ is a transformation of $[0;1]$ preserving
measurability (see [30, 31] for the complex case and [13,14] for the real one;
see also [1], [18], [22], [29]. As for symmetric sequence spaces, see [23,
Ch.IX], [2], [6], [8]). Furthermore, $\phi$ should be measure--preserving
except in the case where $E$ coincides with some of the $L_p$, probably endowed
with a new equivalent norm [30] (see also [14], [18], [22]). In particular,
this shows that $L_p$ are the only r.i. spaces where the orbits of the isometry
group are dense in the unit sphere.\\

The content of the paper is the following. Section 1 contains a preliminary
finite dimensional version of Theorem 1. The proofs of Theorems 1 and 2 are
given in section 2. Besides, section 2 contains a version of Theorem 1 where no
 operator topology is prescribed (see Theorem 2.10). In section 3 we classify
the Coxeter groups in infinite dimensional case (probably, this classification
is not new). In sections 4 and 5 we consider Banach spaces possessing total
families of isometric reflexions. A kind of a structure theorem for isometry
groups is proven (Theorem 5.7). It applies the notions of Hilbert and Coxeter
partial orthogonal subspace decompositions, introduced earlier in this section.
In the last section we give an application to isometry groups of the ideal
generalized sequence spaces.\\

The main results of this paper were announced in [26]; see [27] for their
proofs. Somehow, the proofs have never been published before. The present
article contains some new facts, and the exposition of the old ones is quite
different.

\bigskip

\section{ Isometric reflections in finite dimensional Banach spaces}

Let $A$ be a set of reflexions in a real vector space $E$ and $W$ be the group
generated by the reflexions in $A$. Denote by $\Gamma_{W, A}$ the Coxeter graph
of $W$. Recall [5] that $\Gamma_{W, A}$ has $A$ as the set of vertices; two
vertices are connected by an edge iff the corresponding reflexions do not
commute. By   $\Gamma_W$ we denote the full Coxeter graph of $W$, i.e.  $
\Gamma_W = \Gamma_{W,R}$, where $R = R(W)$ is the set of all the reflexions in
$W$.

\bigskip

\noin {\bf 1.1. Lemma} ([5, Ch. V, 3.7]). {\sl A group $W$ generated by a set
$A$ of orthogonal reflexions in ${\bf R}^n$ is irreducible iff the origin is
the only fixed point of $W$ and the Coxeter graph $\Gamma_{W, A}$ is connected.
In particular, $\Gamma_W$ is connected iff its subgraph $\Gamma_{W, A}$ is
connected.}  \bigskip

Let $E$ be a finite dimensional Banach space. Then ${\rm Iso}\, E$ is a compact
Lie group, and there exists a scalar product in $E$ invariant with respect to
${\rm Iso}\, E$. It can be defined, for instance, by averaging of any given
scalar product over the Haar measure on ${\rm Iso}\, E$. In general, such an
invariant scalar product is not unique. Being orthogonal, two isometric
reflexions in $E$ along the vectors $e_1, e_2 \in S(E)$ commute iff either $e_1
= \pm e_2$ or $e_1 \perp  e_2$. \\

The proof of the following lemma is simple and can be omited.  \bigskip

\noin {\bf 1.2. Lemma.} {\sl Let a connected submanifold $M$ of ${\bf R}^n$ be
invariant under a reflexion $s_{e,e^*}$ which fixes a point $x \in M$ and acts
identically on the tangent space $T_x M$. Then $M$ is contained in the mirror
hyperplane ${\rm Ker}\, e^*$. }

\bigskip

The main result of this section is the following

\bigskip

\noin {\bf 1.3. Proposition.} {\sl Let $E$ be a real Banach space of dimension
$n$. Let $G \subset {\rm Iso}\, E$ be a closed subgroup of a positive dimension
which contains reflexions $t_1,\dots,t_n$ along linearly independent vectors
$e_1, \dots,e_n$. Then there exists a subspace $H \subset E$ such that

\noin a) ${\rm dim} H \ge 2, \,\,\, H$ is euclidean and biorthogonally
complemented in $E$;

\noin b) the unit sphere $S(H)$ of $H$ coincides with an orbit of the identity
component $G_0$ of $G$;

\noin c) there exists a projection $p$ of $E$ onto $H$ such that $1_E - 2p \in
G$ and $p$ commutes with any reflexion $t \in G$. Furthermore, $(1_E - p) +
{\bar u}p \in G$ for any ${\bar u} \in O(H)$.} \bigskip

\noin {\sl Proof.} Fix an invariant scalar product in $E$ and identify $E$ with
${\bf R}^n$ in such a way that ${\rm Iso}\, E \subset O(n)$. Let $u_1, \dots,
u_n$ be the system of vectors in ${\bf R}^n$ biorthogonal to the system $e_1,
\dots, e_n$. Since $\dim G > 0$, the orbit $Gu_i$ has a positive dimension for
at least one value of $i$, say for $i=1$. We may also assume that $u_1 \in
S^{n-1}$, where $S^{n-1}$ is the euclidean unit sphere in ${\bf R}^n$. Let $M$
be the connected component of the orbit $G u_1$ which contains $u_1$. Since
$u_1$ is fixed by any of the reflexions $t_i, i= 2,\dots, n$, $M$ is invariant
under these reflexions, an hence the tangent space $T = T_{u_1} M$ is
invariant, too. Thus for each $i = 2,\dots, n$ either $e_i \in T$ or $e_i \perp
T$. Put $$A = \{i \in \{2, \dots, n\}\, | \,e_i \in T\}$$ and $$B = \{i \in
\{2, \dots, n\} \,|\, e_i \perp T\}\,\,\, .$$  Since $G \subset O(n)$ and $u_1
\in S^{n-1}$, we have $M \subset S^{n-1}$, and so $T \subset T_{u_1} S^{n-1}$.
Therefore $T\perp  u_1$. It follows that $T \subset L$, where $L = {\rm span}
(e_2, \dots, e_n)$, and therefore $T = {\rm span}(e_i\,|\, i \in A)$
(hereafter {\it span} means {\it the linear span}).

Thus, $\dim M = \dim T = {\rm card} A$. Since $M$ is $t_i$-invariant for $i \in
B$, by Lemma 1.2, $M$ is contained in the subspace $H = \{v \in E \,|\, v \perp
e_i, i \in B\}$. It is easily seen that $k= \dim H = {\rm card} A + 1 = \dim M
+ 1 $. Thus $M$ is a closed submanifold of each of the unit spheres $S_r (H) =
S_r (E) \cap H$, where $r= ||u_1 ||_E$, and $S^{k-1} = S^{n-1} \cap H$, of the
same dimension $\dim M = \dim H - 1 = k - 1$. Hence $M$ coincides with both of
them. At the same time, being connected $M$ coincides with the orbit $G_0 u_1
$. Here $k \ge 2$, since $\dim M > 0$. Therefore, $H$ is euclidean and the unit
sphere $S(H)$ coincides with the orbit $G_0 (u_1 /r)$.

Since $T \subset H$ and $e_i \in T$ for each $i \in A$, where ${\rm card} A =
k-1>0$, there exists $i_0 \in A$ such that $e_{i_0} \in S(H)$, and thus $S(H) =
G_0 e_{i_0}$. Let $w_1, \dots , w_k \in H$ be an orthogonal basis in $H$ with
$||w_i||_E =1, i=1, \dots, k$, and $g_1, \dots, g_k \in G$ be such that $g_j
(e_{i_0}) = w_j$. Then $s_j = g_j t_{i_0} g_j ^{-1} \in G$ is the orthogonal
reflexion along the vector $w_j, j=1, \dots, k$. By the same reasoning as
above, for any vector $w \in S(H)$ the orthogonal reflexion $s_{w, w^*}$ along
$w$ belongs to $G$.

The reflexions $s_j , j=1, \dots , k$, pairwise commute, and so $p= {1\over
2}(1_E -   \prod _{i=1}^k s_i)$ is the orthogonal projection of $E$ onto $H$
such that $\tau =1_E - 2p = \prod  _{i=1}^k s_i \in G \subset {\rm Iso}\, E$.
Thus, $||p||_E = {1 \over 2}||1_E - \tau||_E =1$ and either $E=H$ or $||1_E -
p||_E = {1 \over 2} ||1_E + \tau||_E = 1$. Therefore, $H$ is a biorthogonally
complemented subspace of $E$.

Any orthogonal reflexion $\bar s $ in $H$ coincides with the restriction to $H$
of some reflexion $s \in G$, where in fact $s = (1_E - p) + {\bar s}p$. The
same is true for any orthogonal operator ${\bar u} \in O(H)$; indeed, the group
$O(H)$ is generated by orthogonal reflexions.

Let $t \in G$ be a reflexion. The mirror hyperplane of $t$ intersects with $H$
by a subspace of $H$ of dimension $k-1 >0$. Therefore, $t$ has a fixed point on
the sphere $M=S^{k-1} \subset H$, and so $t(M) \cup M$ is connected and
contained in the orbit $Gu_1$. It follows that $t(M) = M$, $H$ is invariant
with respect to $t$ and so $t$ and $p$ commute. This completes the proof.
\,\,\,$\bigcirc$

\bigskip

\noin {\bf 1.4. Corollary.} {\sl Let $W$ be a group generated by isometric
reflexions in a finite dimensional Banach space $E$. If $W$ is irreducible and
infinite, then $E$ is euclidean and $W$ is dense in the orthogonal group ${\rm
Iso}\, E \approx {\it O}(n), n=\dim E$. } \bigskip

\noin {\sl Proof.} Let $G$ be the closure of $W$ in ${\rm Iso}\, E$ and $G_0$
be the identity component of $G$. Sinse $W$ is irreducible, by Lemma 1.1, it
contains $n$ reflexions along linearly independent vectors, and the Coxeter
graph $\Gamma _W$ is connected. Let $H$ be the euclidean subspace of $E$
constructed in  Proposition 1.3. Since by (c), $H$ is invariant with respect to
the reflexions from $W$, for each $s_{e,e^*} \in W$ either $e \in H$ or $e
\perp H$. If $A$ resp. $B$ is the set of reflexions from $W$ of the first resp.
second type, then each element of $A$ commutes with every element of $B$. By
the connectedness of the graph $\Gamma _W$ one of the sets $A$ and $B$ should
be empty. This shows that $H=E$. By (c), ${\bar u} \in G$ for any ${\bar u} \in
O(H)$. Therefore, $G = O(H)$ and we are done. \,\,\, $\bigcirc$

\bigskip

\noin {\it Remark.} Related results can be found in [6], [10, (1.7)], [23],
[25].

\section{ Proofs of Theorems 1 and 2}

\smallskip

\noin {\bf 2.1. Definition.} Let E be a real Banach space, and let $s_1, s_2$
be two isometric reflexions in $E$ along linearly independent vectors $e_1$,
$e_2 \in S=S(E)$.  Denote by $\alpha (s_1 , s_2)$ the minimal positive angle
between the lines containing $e_1$ and $e_2$, measured with respect to an
invariant inner product in the plane $L={\rm span}(e_1, e_2)$. Put $\alpha (s_1
, s_2) = 0$ iff $e_1 = \pm e_2$.  \bigskip

\noin {\bf 2.2.} {\sl Remarks}.
a) It is easily seen that the above definition does not depend on the choice of
an invariant scalar product in $L$.

\noin b) An isometric reflection $s=s_{e,e^*}$ in $E$ is uniquely defined by
the reflexion point $e \in S(E)$. Indeed, this is true for the restriction of
$s$ to any finite dimensional subspace $F$ containing $e$, since the mirror
hyperplane ${\rm Ker}\, e^* \cap F$ of $s|F$ is orthogonal to $e$ with respect
to an invariant scalar product on $F$.  Thus, this is true for $s$ itself.

\noin c) Two isometric reflections $s_1$ and $s_2$ commute iff either $\alpha
(s_1 , s_2) = 0$, i.e. $e_1 =  \pm e_2$ or $\alpha (s_1 , s_2) = {\pi \over
2}$, i.e. $e_1  \perp  e_2$ in $L$.

\bigskip

\noin {\bf 2.3. Lemma.} {\sl Let $s_i = s_{e_i ,e{_i}^*}, i=1,2$, be two
isometric reflexions in $E$. Then} $$\cos ^2 \alpha (s_1 , s_2) = {e_1}^*
(e_2){e_2}^* (e_1) .$$
{\sl Proof. }This is evidently true if $e_1 = \pm  e_2$. Assume, further, that
$e_1$ and $e_2$ are linearly independent. Let an invariant scalar product in
the plane $L={\rm span}(e_1 , e_2)$ be given by the bilinear form $B = \left(
\begin{array}{cc} b & a\\ a & c \end{array} \right)$ with respect to the basis
$(e_1 , e_2)$ in $L$. Consider the orthogonal projection $p_i = {1 \over 2}
(1_L + s_i |L)$ of $L$ onto the mirror line $l_i$ of the axial reflexion $s_i
|L \,,i=1,2$. Since $p_i (e_j) \perp  e_i$ for $j \ne i$, we have $$0=B(p_1
(e_2), e_1) = B(e_2 - {e_1}^* (e_2) e_1, e_1) = a - {e_1}^* (e_2)b$$ and $$
0=B(p_2 (e_1), e_2) = a - {e_2}^* (e_1)c\,\, .$$ Thus $$ a^2 = {e_1}^*
(e_2){e_2}^* (e_1)bc \,\,,$$ and so $$\cos ^2 \alpha (s_1 , s_2) = {a^2 \over
bc} = {e_1}^* (e_2){e_2}^* (e_1)\,\, .$$ \,\,\,$\bigcirc$

\bigskip

\noin {\bf 2.4. Corollary.} {\sl  $$\cos \alpha (s_1 , s_2) \ge 1 - ||e_1 -
e_2||_E .$$ In particular, if $e_1 \ne e_2$ and $||e_1 - e_2||_E < 1$, then
$s_1$ and $s_2$ do not commute. } \bigskip

\noin {\sl Proof. } Since $s_i \in {\rm Iso}\, E$, and so $||e_i||_E =
||e^*_i||_{E^*}=e^*(e)=1$, we have $$|1 - e_1^*(e_2)| = |e_1^*(e_1-e_2)| \le
||e_1-e_2||_E$$ and $$|1-e_2^*(e_1)| \le ||e_1-e_2||_E\,\,.$$ We can assume
that $||e_1-e_2||_E < 1$. Then from the above inequalities we obtain
$$|e_1^*(e_2)| \ge 1 - ||e_1-e_2||_E$$ and $$|e_2^*(e_1)| \ge 1 -
||e_1-e_2||_E\,\,.$$ The desired inequality follows from the latter two by
multiplying them and making use of Lemma 2.3. \,\,\,$\bigcirc$

\bigskip

\noin {\bf 2.5. Lemma.} {\sl Let $s=s_{e, e^*} \in {\rm Iso}\, E$. Consider the
function on ${\rm Iso}\, E \times {\rm Iso}\, E$ $$\phi_s (g_1, g_2) = \sin^2
\alpha (s_1 , s_2)\, $$ where $s_i = g_i s {g_i}^{-1},\, i=1,2$. Then

\noin a) $\phi_s$ is left invariant, i.e. $$ \phi_s (g_1, g_2) = \phi_s (gg_1,
gg_2) = \phi_s (1_E , {g_1}^{-1}g_2) $$ for each $g, g_1, g_2 \in {\rm Iso}\,
E$.

\noin b) $$\phi_s (g_1, g_2) = 1 - {e^*}({g_1}^{-1}g_2 (e)) e^* ({g_2}^{-1}g_1
(e)) \,\,.$$ Therefore, $\phi _s$ is continuous on $({\rm Iso}\, E)^2$ in the
strong operator topology.

\noin c) For any two elements $g', g'' \in G_0$ such that $\phi_s (g', g'' ) >
0$, and for any $\epsilon,\, 0< \epsilon < 1 ,$ one can find a finite chain of
elements $h_0 = g', h_1, \dots, h_n = g''$ with the property $0<\phi_s (h_i ,
h_{i+1}) < \epsilon$, so that the reflexions $t_i = h_i s {h_i}^{-1}$ and
$t_{i+1}$ do not commute for all $i = 0, 1, \dots, n-1$.}\bigskip

\noin \noin {\sl Proof}. (a) is evident. The identity in (b) easily follows
from the equality $$\phi_s (g_1, g_2) = 1 - {e^*}({g_1}^{-1}g_2
(e))({g_1}^{-1}g_2)^* (e^* )(e)\,\,\, ,$$ which follows from (a) and Lemma 2.3.
The second statement of (b) is true since ${\rm Iso}\, E$ is a topological
group with respect to the strong operator topology. To prove (c), consider the
covering of $G_0$ by the open subsets $$ U_{\epsilon} (g) = \{h \in G_0 \,
\vert \, \phi_s (g, h) < \epsilon\}\,\,\, .$$ Since $G_0$ is connected, any two
of them $ U_{\epsilon} (g')$ and $ U_{\epsilon} (g'')$ can be connected by a
finite chain of such subsets, and the assertion follows.
\,\,\,$\bigcirc$\bigskip

\noin {\bf 2.6. Proposition. }{\sl Let $s=s_{e, e^*} \in {\rm Iso}\, E, \,\,
g_1, \dots, g_n \in G_0$ and $H' = {\rm span} (e_1, \dots, e_n)$, where $e_i =
g_i (e),\, i=1,\dots, n$. Then

\noin a) $H'$ is euclidean;

\noin b) there exists a unique projection $p'$ of $E$ onto $H'$ such that $1_E
- 2p \in {\rm Iso}\, E$;

\noin c) the unit sphere $S(H')$ of $H'$ is contained in the orbit $G_0\, e$,
and for each vector $v \in S(H')$ there exists a reflexion $s_{v, v^*} \in {\rm
Iso}\, E$  along $v$ commuting with $p'$.} \bigskip

\noin {\sl Proof.} First we construct a finite dimensional subspace $F$
containing $H'$ which satisfies all the properties of (a), (b), (c) above.

Put $g_0 = 1_E$ and for each pair $(g_i , g_{i+1}), i=0, \dots, n-1$, find a
chain $\{h_{ij}\}_{j=0}^{n_i}$ as in Lemma 2.5.c above. The proposition is
evident in the case when $\dim H' =1$ , and so we may assume that $g_i (e) \ne
e$ for at least one value $i_0$ of $i$. Since the continuous function $\phi_s$
takes all its intermediate values on $G_0$, we can also choose the element
$h=h_{i_0 , 1}$ in such a way that the angle $\alpha (s, hsh^{-1})$ is
irrational modulo $\pi$, and thus the group generated by the reflections $s$
and $hsh^{-1}$ is infinite.

Put $F = {\rm span}(h_{ij} (e)\,\vert \, j=0,\dots, n_i,  i=0,\dots,n)$. Let
$W$ be the group generated by the reflexions $\{t_{ij}|F\}$ in $F$, where
$t_{ij} = h_{ij} s h_{ij}^{-1}, j=0,\dots, n_i,  i=0,\dots,n$. It is clear that
the origin is the only fixed point of $W$ in $F$. Since, by the construction
the Coxeter graph, $\Gamma_W$ of $W$ is connected, by Lemma 1.1, $W$ is
irreducible. $W$ being infinite, by Corrolary 1.4, the subspace $F$ is
euclidean and the closure of $W$ coincides with the group ${\rm Iso}\, E =
O(F)$.

Let $v_1,\dots,v_l$, where $l=\dim F$, be a basis of $F$ chosen from the system
$(h_{ij} (e))$, and $t_k = t_{v_k , {v_k}^*}, k=1,\dots,l$, be the
corresponding reflexions from the system $(t_{ij})$. Put $M= \bigcap _{k=1} ^l
{\rm Ker}\, {v_k}^*$. It is easily seen that $E = M \oplus F$. Let $F'$ be a
finite dimensional subspace of $E$ containing $F$, endowed with an invariant
scalar product. Then for each $k=1,\dots,l$ the restriction $t'_k = t_k |F'$ is
an orthogonal reflexion in $F'$, and so $v_k \perp  ({\rm Ker} \, {v_k}^* \cap
F')$. Therefore, $F \perp (M \cap F')$. It follows that each of the vectors
$h_{ij} (e) \in F$ is orthogonal to $M \cap F'$, too, so that the restriction
$t_{ij} |(M \cap F')$ is the identity mapping. This gives the representation
$t_{ij} = (1_E - p) + t_{ij} p$, where $p$ is the projection of $E$ onto $F$
along $M$. Thus, each element $ \bar g \in W$ can be represented as the
restriction to $F$ of the isometry $g = (1_E - p) + {\bar g}p \in {\rm Iso}\,
E$. If a sequence ${\bar g}_i \in W$ converges to an element ${\bar h} \in
O(F)$, then the sequence of extensions $g_i$ converges to the extension $h =
(1_E - p) + {\bar h} p$ of ${\bar h}$, where $h \in\, {\rm Iso}\, E$. In
particular, in this way each orthogonal reflexion in $F$ extends to a unique
isometric reflexion in $E$, and each element ${\bar u} \in O(F)$ extends to the
unique isometry $u=(1_E - p) + {\bar u}p \in {\rm Iso}\, E$. It follows that
$S(F) \subset G_0 (e)$.

Let $f_1,\dots, f_l$ be an orthogonal basis in $F$ and ${\bar s}_1, \dots,
{\bar s}_l$ be the orthogonal reflexions in $F$ along these vectors. It is
easily seen that then $p = {1 \over 2} (1_E - \prod _{i=1}^l s_i)$, and thus
$1_E - 2p = \prod _{i=1}^l s_i \in {\rm Iso}\, E$. If $s' \in {\rm Iso}\, E$ is
a reflexion along a vector $v' \in S(F)$, then as above $s' = (1_E - p) + s'p =
(1_E - p) + ps'p$, and so $s'$ and $p$ commute.

It is evident that the subspace $H' \subset F$ has the same properties as $F$
itself, and therefore (a), (b), (c) are fulfilled.\,\,\, $\bigcirc$

\bigskip

\noin {\bf 2.7.} {\sl Remark.} It is easily seen that if $H' \subset H''$ are
two subspaces as in Proposition 2.6, then for the corresponding projections $p'
, p''$ we have $p' \prec p''$, i.e. $p'p'' = p' (= p''p')$. \bigskip

\noin {\bf 2.8.} {\sl Proof of Theorem 1.a}.  Let $x, y$ be two arbitrary
vectors in $H$. Then for any $\epsilon > 0$ in the linear span of the orbit
$G_0 e$ there exist two vectors $x_{\epsilon}, y_{\epsilon}$ such that $||x -
x_{\epsilon}||_E < \epsilon, ||y -  y_{\epsilon}||_E < \epsilon$. Let
$x_{\epsilon} = \sum _{i=1}^n a_i g_i (e)$ and  $y_{\epsilon} = \sum_{i=1}^n
b_i g_i (e)$, where $g_i \in G_0, i = 1, \dots, n$. Put $H' = {\rm span} (g_i
(e)\vert i = 1, \dots, n)$. By Proposition 2.6, the  subspace $H'$ is
euclidean, and therefore the norm in $H'$ satisfies the four squares identity.
In particular, $$||x_{\epsilon} + y_{\epsilon}||^2 + ||x_{\epsilon} -
y_{\epsilon}||^2 = 2 (||x_{\epsilon}||^2 + ||y_{\epsilon}||^2).$$ Passing to
the limit we see that the same identity holds for $x, y \in H$. It follows that
$H$ is a Hilbert space (see [7], Ch.7, \S 3).

 Consider further the family of all finite dimensional subspaces ${H'}$ which
belong to the linear span of the orbit $G_0 e$. Let ${\cal P} = \{ p' \}$ be
the corresponding partially ordered family of finite dimensional projections $E
\rightarrow H'$ such that $1_E - 2p' \in {\rm Iso}\, E$. For a fixed vector $v
\in S(E)$ and for each $p' \in P$ consider the subset $$Y_{p'} = \omega
\{p''(v)\,\, |\,\, p'' \in P, \,\, p' \prec p''\}\,\, ,$$ where $\omega$
denotes the closure with respect to the weak topology in $E$. The family
$\{Y_{p'}\}$ has the property that for each finite system of projections $p'_1,
\dots, p'_n \in P$ the intersection $\cap _{i=1} ^n Y_{p'_i}$ is non-empty.
Indeed, let $H_0 = {\rm span} ({\rm Im}\, p'_1, \dots, {\rm Im}\, p'_n)$, and
let $p'_0 \in \cal P$ be the corresponding projection of $E$ onto $H_0$. Then
$p'_i \prec p'_0$, and hence $p'_0(v) \in Y_{p'_i}$ for each $i=1, \dots, n$.
Since $H$ is a Hilbert space, the unit ball $B(H)$ of $H$ is weakly compact. It
follows that the centralized family $\{Y_{p'}\}_{p' \in \cal P}$ of weakly
closed subsets of $B(H)$ has a non-empty intersection. By the Barry's Theorem
[4] the generalized sequence of projections ${\cal P} = (p')$ converges in the
strong operator topology to its upper bound $p$ which is a projection of $E$
onto $H$, and which satisfies the condition ($i$) $1_E - 2p \in {\rm Iso}\, E$.
In particular, $H$ is biorthogonally complemented in $E$. This proves ($a$).

\noin $b, \,\,c$. Let $s' = s_{x, x^*} \in {\rm Iso}\, E$, where $v \in G_0 e$.
Then $s'$ commutes with any projection $p' \in \cal P$ such that $p'(x) = x$,
which means that ${\rm Ker}\, p' \subset {\rm Ker}\, x^*$. Passing to the limit
we see that $p$ commutes with $s'$ and ${\rm Ker}\, p \subset {\rm Ker}\, x^*$,
too. It follows that $s' = (1_E - p) + s'p$.

Let $x_0 \in S(H)$ be the limit of a generalized sequence of vectors
$x_{\alpha} \in G_0 e \cap S(H)$. Then the corresponding sequence of isometric
reflections $s_{\alpha} = s_{x_{\alpha}, x^* _{\alpha}} = g_{\alpha} s
g_{\alpha}^{-1}$, where $g_{\alpha} \in G_0$ and $g_{\alpha} (e) = x_{\alpha}$,
is strongly convergent to the reflexion $s_0 = s_{x_0 , x^* _0} \in {\rm Iso}\,
E$. Indeed, from the representation $s_{\alpha} = (1_E - p) + s_{\alpha} p$ it
easily follows that the generalized sequence $\{s_{\alpha}\}$ converges to $s_0
= (1_E - p) + s_0 p$ on each of the complementary subspaces ${\rm Ker}\, p$ and
$H$.

Let ${\bar u} \in O(H)$ and $x \in E$. Consider the extension $u = (1_E - p) +
{\bar u}p$ of ${\bar u}$ to $E$. Since $||p(x)||_E = ||up(x)||_E$, there exists
an orthogonal reflexion ${\bar s}_0$ in $H$ such that ${\bar s}_0 p(x) =
up(x)$. Let $s_0 = (1_E -p) + {\bar s}_0 p \in {\rm Iso}\, E$. Then we have
$||x||_E = ||s_0 (x)||_E = ||(1_E - p) (x) + {\bar s}_0 p (x)||_E =
||u(x)||_E$. Therefore, $u \in {\rm Iso}\, E$, and thus ($ii$) is fulfilled.

Now it is clear that the orbit $G_0 e \subset S(H)$ contains the orbit of the
strong identity component of the orthogonal group $O(H)$, and so it coincides
with the sphere $S(H)$. This proves ($c$).

Let $g \in {\rm Iso}\, E$ be such that ${\bar u} = g|H \in O(H)$. We will show
that $g$ leaves the subspace ${\rm Ker}\, p$ invariant and thus ${\bar v} = g|
{\rm Ker}\, p \in {\rm Iso} ({\rm Ker}\, p)$ and $g = {\bar v} (1_E - p) +
{\bar u} p$.

Suppose that $g({\rm Ker}\, p)   \not\subset    {\rm Ker}\, p$. Consider the
operator $g_1 = g u^{-1} \in {\rm Iso}\, E$. We have $g_1|H = 1_H$ and
$g_1|{\rm Ker}\, p = g|{\rm Ker}\, p$. By our assumption there exists a vector
$x \in {\rm Ker}\, p$ such that $g_1 (x)  \notin    {\rm Ker}\, p$, and so
$pg_1 (x) \ne 0$. Denote $x_1 = (1_E - p)g_1(x)$ and $x_2 = pg_1(x)$. Then
$g_1(x) = x_1 + x_2$, hence ${g_1}^{-1} (x_1) = {g_1}^{-1}(g_1(x) - x_2) = x -
x_2$.

Consider two functions $\phi (t) = ||x_1 + tx_2||_E$ and $\psi (t) = ||x +
tx_2||_E$. Since $1_E - 2p \in {\rm Iso}\, E, x, x_1 \in {\rm Ker}\, p$ and
$x_2 \in {\rm Im}\, p$, we have $||x_1 + x_2|| = ||x_1 - tx_2||$ and $||x +
tx_2|| = ||x - tx_2||$. Thus, both $\phi$ and $\psi$ are even functions. From
the equalities $${g_1}^{-1} (x_1 + tx_2) = x - x_1 + tx_2 = x + (1 - t)x_2$$
and $$g_1 (x + tx_2) = x_1 + (1 + t)x_2$$ and the fact that $g_1 \in {\rm
Iso}\, E$ we obtain that $\phi (t) =  \psi (1 - t)$ and $\psi (t) =  \phi (1 +
t)$. It follows that $ \phi(t) = \phi (-t) = \psi (1 + t) = \phi (t + 2)$.
Therefore,  being convex and periodic function on $\bf R$, $\,\,\phi (t)$
should be constant. This is possible only if $x_2 = 0$, i.e. $g_1 (x) \in {\rm
Ker}\, p$, which is a contradiction.  Thus, ($iii$) is fulfilled as well. This
completes the proof of Theorem 1.\,\,\, $\bigcirc$

\bigskip

It has been already noted that the statement of Theorem 2.b is a direct
corollary of Theorem 1.  Thus, it is enough to prove Theorem 2.a.

\bigskip

\noin {\bf 2.9.} {\sl Proof of Theorem 2.a}. It is enough to show that the four
squares identity holds in E. For this it is enough, as it was done in the proof
of Theorem 1, to approximate an arbitrary pair of vectors $x, y \in E$ by a
sequence $\{(x_{\alpha} , y_{\alpha})\}$ of pairs of vectors belonging to
finite dimensional euclidean subspaces $H_{\alpha}$ of $E$. In turn, it is
enough to show that any pair of vectors in the linear span of the orbit $Ge$ of
the group $G = {\rm Iso}\, E$ belongs to a finite dimensional euclidean
subspace $H'$ of $E$. Indeed, it is easily seen that under our assumptions any
orbit of $G$ in the unit sphere $S = S(E)$ is dense in $S$. In particular, the
orbit $Ge$ is dense in $S$.

Fix such a pair $x, y \in {\rm span}(Ge)$ and consider a subspace
$H' = {\rm span} (g_1(e),\dots, g_n(e))$, $g_i \in G, i=1,\dots,n,$
containing this pair. Since the orbit $Ge$ is dense in $S$, for any two vectors
$g'(e)$ and $g''(e) \in G\,e$ one can find a finite chain of vectors
$h_j (e) \in G\,e, j=0,\dots, k, $
such that
$h_0 (e) = g'(e), h_k (e) = g''(e)$
and
$||h_{j + 1} (e) - h_j (e)||_E < {1 \over 10},\,\, j=0,1,\dots,k-1$.
 Find such a chain
$\{h_{ij} (e)\} _{j=1,\dots,k_i}$
 for each of the pairs
 $(g_i (e) , g_{i+1} (e)), i=1,\dots,n-1,$
and put
 $F={\rm span}(h_{ij} (e), j=0,\dots,k_i , i=1,\dots,k-1)$.
 We may assume that $E$ is infinite dimensional (otherwise the proof is
simple), and that $\dim L > 8$. Let
$\{t_{ij} = h_{ij} s {h_{ij}}^{-1}\}$
be the system of isometric reflexions along the vectors
$h_{ij} (e) , j=1,\dots,k_i , i=1,\dots,n-1$,
and let $W$ be the group generated by the restrictions $t_{ij}\vert F$. By
Corollary 2.4, the Coxeter graph ${\Gamma}_W$ is connected, and since the
system of vectors $(h_{ij} (e))$ is complete in $F$, by Lemma 1.1, the group
$W$ is irreducible. For a pair of vectors $(v' = h_{ij} (e), v'' = h_{ij+1}
(e))$ we have $0 < ||v' - v''|| < {1 \over 10}$, and so by Corollary 2.4, $$0 <
{\alpha} (t_{ij} , t_{ij+1} ) < \arccos {9 \over 10}\, .$$ From the
classification of Coxeter groups [5, Ch.VI, sect.4] it follows that the group
$W$ is infinite. Thus, by Corollary 1.4, the subspaces $F$ and $H' \subset F$
are euclidean. The theorem is proven. \,\,\,$\bigcirc$\\

\bigskip

{\it A priori}, the strong operator topology could be still too strong in order
that the identity component $G_0$ be big enough to apply Theorem 1 in an
efficient way.  Next we give a version of Theorem 1 which does not involve any
operator topology.

Recall that a group $G$ is {\it locally finite} if every finitely generated
subgroup of $G$ is finite.

\bigskip

\noin {\bf 2.10. Theorem.} {\it Let $s = s_{e, e^*}$ be an isometric reflexion
in a Banach space $E$. Denote $$U = \{g \in {\rm Iso}\,E\,|\,[s, g^{-1}sg] \neq
1_E\}\,\,.$$ Let $G_1$ be the subgroup of ${\rm Iso}\,E$ generated by $U$, and
$H = {\rm \overline {span}}\,(G_1 \, e)$. If any of the following two
conditions (i), (ii) is fulfilled, then all the conclusions (a), (b), (c) of
Theorem 1 hold:\\
i) The group $W_1$ generated by the set of reflexions ${\rm IR}_1 =
\{g^{-1}sg\}_{g \in G_1}$ is not locally finite.\\
ii) The orbit $G_1 \,e$ contains three linearly independent vectors $e_1 , e_2
, e_3$, where $||e_1 - e_2 ||_E < 1 - \cos \pi/5$.}

\bigskip

\noin \noin {\sl Proof}. Repeating the arguments used in the proofs of Theorems
1 and 2, it is enough to show that for each finite subset $\sigma \subset {\rm
IR}_1$ there exists another finite subset $\sigma_1 \subset {\rm IR}_1$ such
that $\sigma \subset \sigma_1$ and the group generated by the reflexions from
$\sigma_1$, as well as its restiction to the subspace ${\rm span}\,(v\,|\,s_{v,
v^*} \in {\rm IR}_1)$, is infinite. In other words, each finite subgraph
$\gamma$ of the Coxeter graph $\Gamma_{W_1}$ should be contained in a finite
connected subgraph $\gamma_1 \subset \Gamma_{W_1}$ with the following property:
the group $W(\gamma_1)$ generated by the reflexions which correspond to the
vertices of $\gamma_1$, is infinite. The latter holds as soon as the Coxeter
graph $\Gamma_{W_1}$ is connected and contains a finite subgraph $\gamma_0$
such that the group $W(\gamma_0)$ is infinite. If the first of these conditions
is fulfilled, than the second follows from each of the assumptions  ($i$) and
($ii$) above. Indeed, it is clear for ($i$). As for ($ii$), by the
connectedness of the graph $\Gamma_{W_1}$, one can find a finite connected
subgraph $\gamma_0 \subset \Gamma_{W_1}$ which contains three vertices
corresponding to the reflexions $s_1 , s_2 , s_3 \in {\rm IR}_1$ with the
reflexion vectors $e_1 , e_2 , e_3$, resp. From the classification of finite
Coxeter groups [5, Ch.VI, sect.4] it follows that the group $W(\gamma_0)$ is
infinite. Indeed, if $V(\gamma_0)$ is the subspace generated by the reflexion
vectors of the reflexions from $W(\gamma_0)$, then ${\rm dim}\,V(\gamma_0) \ge
3$ and by Corollary 2.3, the order of the rotation $s_1 s_2 \in W(\gamma_0)$ is
greater than $5$.

 Thus, it remains to check that the graph $\Gamma_{W_1}$ is connected. Let the
vertices $v, v'$ of $\Gamma_{W_1}$ correspond to the reflexions $s, s' = h^{-1}
s h$ resp., where $h = g_n g_{n-1} \cdot\dots\cdot g_1 \in G_1$ is arbitrary
and $g_i \in U\,,\,i=1,\dots,n$. Put $h_i = g_i g_{i-1}\cdot\dots\cdot g_1$ and
$s_i = h_i^{-1} s h_i\,,\,i=0,\dots,n$, so that $h_0 = 1_E , h_1 = g_1 , h_n =
h$ and $s_0 = s\,,\,s_n = s'$. Since $g_1 \in U$, the reflexions $s_0 = s$ and
$s_1$ do not commute, and thus $0 < \phi (1_E , h_1 ) < 1$. By Lemma 2.4.a, $0
< \phi (1_E , h_1 ) = \phi (g_1 , g_1 h_1 ) = \phi (h_1 , h_2) < 1$ , and
therefore the reflexions $s_1$ and $s_2$ do not commute, as well. By induction,
we see that $s_i$ does not commute with $s_{i+1}$ for all $i=0,\dots, n-1$, and
so the vertices $v, v'$ of the graph $\Gamma_{W_1}$ are connected by a path.
This concludes the proof. \,\,\,$\bigcirc$

\section{ Infinite Coxeter groups}

In sections 4, 5, 6 below we will use the classification of infinite Coxeter
groups. Although it should be well known, in view of the lack of references we
reproduce it here in all details.

  By {\it an infinite Coxeter group} we mean an infinite locally finite group
$W$ generated by reflexions in a real vector space $V$ which is algebraically
irreducible in $V$. We fix the following notation and conventions.\\

\noin {\bf 3.1.} {\sl Notation}. Denote by ${\bf R}^{\Delta}$ the linear space
of all the real functions with finite support defined on a given set $\Delta$,
and by ${\bf R}^{\Delta}_0$ the subspace of functions with the zero mean value.
Let

\smallskip

\noin $A_{\Delta}$ be the group of finite permutations of elements of $\Delta$
acting in ${\bf R}^{\Delta}_0$;

\smallskip

\noin $B_{\Delta}$  be the group of finite permutations of $\Delta$ and changes
of sign of values at the points of finite subsets of $\Delta$ acting in ${\bf
R}^{\Delta}$;

\smallskip

\noin $D_{\Delta}$ be the subgroup of $B_{\Delta}$ which consists of finite
permutations and changes of signs of even numbers of coordinates
acting in ${\bf R}^{\Delta}$.

\smallskip

If $\Delta$ is infinite then $A_{\Delta}$, $B_{\Delta}$, $D_{\Delta}$ are
infinite Coxeter groups. Let $\epsilon_{\delta}$ be the characteristic function
of the one-point subset $\{\delta\}$ of $\Delta$, so that
$(\epsilon_{\delta}\,\vert\,\delta \in \Delta)$ is the standard Hamel basis of
${\bf R}^{\Delta}$. Let ${\bf R}^{\Delta}$ be endowed with the standard scalar
product. Then $A_{\Delta}$  (resp. $B_{\Delta}$, $D_{\Delta}$) is generated by
orthogonal reflexions along vectors from the infinite root system
$(\epsilon_{\delta} - \epsilon_{\delta'})$ (resp. $(\pm \epsilon_{\delta},
\,\,\,\pm \epsilon_{\delta} \pm \epsilon_{\delta'})$, $(\pm \epsilon_{\delta}
\pm \epsilon_{\delta'})$) \,\,\,($\delta, \delta' \in \Delta, \,\,\delta \ne
\delta'$).

In the category of pairs $(W, V)$, where  $W$ is a group generated by
reflexions in a real vector space $V$, there is a natural notion of
isomorphism. We will also use a notion of subpair. Namely, we will say that
$(W', V')$ is {\it a subpair} of $(W, V)$ if $W'$ is the restriction of a
subgroup of $W$ generated by reflexions to its invariant subspace $V'$. An
embedding of pairs is an isomorphism with a subpair. In the proposition below
{\it isomorphism of Coxeter groups} means isomorphism of pairs, rather then
isomorphism of abstract groups.

\bigskip

\noin {\bf 3.2. Proposition.} {\sl Any infinite Coxeter group $W$ is isomorphic
to one and only one of the groups $A_{\Delta}$, $B_{\Delta}$, $D_{\Delta}$.}

\bigskip

\noin {\sl Proof}. In the sequal $\gamma$ denotes a finite connected subgraph
of the Coxeter graph $\Gamma_W$, $G(\gamma)$ denotes the finite subgroup of $W$
generated by reflexions $s_i = s_{e_i, e^*_i} \in \gamma, i=1,\dots, {\rm
card}\,\gamma$, $V(\gamma) = {\rm span}\,(e_i\,,i=1,\dots, {\rm card}\,\gamma)$
and $n(\gamma) = \dim\,V(\gamma)$. By Lemma 1.1, the restriction
$G(\gamma)\,\vert\,V(\gamma)$ is irreducible, so it is a finite Coxeter group.
The full Coxeter graph $\Gamma_{G(\gamma)}$ can be naturally identified with a
finite connected subgraph $\bar \gamma$ of $\Gamma_W$ containing $\gamma$; in
fact, $\bar \gamma$ is the maximal subgraph of $\Gamma_W$ with the properties
that $V(\bar \gamma) = V(\gamma)$ and $G(\bar \gamma) = G(\gamma)$ (but the
first one alone does not determine $\bar \gamma$). If $n = n(\gamma)>8$, then
$G(\gamma)$ is one of the Coxeter groups $A_n, B_n, D_n$ [5, Ch.VI, sect. 4].

Let $\Delta$ be a set with ${\rm card}\,\Delta = {\rm dim}\,V$, where ${\rm
dim}\,V$ is the cardinality of a Hamel basis in $V$. The proposition follows
from the assertions (i) - (iii) below.

\smallskip

\noin i) $(W, V) \approx (B_{\Delta}, {\bf R}^{\Delta})$ if $(G(\gamma_0),
V(\gamma_0))  \approx (B_n, {\bf R}^n)$ for some $\gamma_0 \subset \Gamma_W$;

\noin ii) $(W, V) \approx (D_{\Delta}, {\bf R}^{\Delta})$ if there is no
$\gamma \subset \Gamma_W$ for which $(G(\gamma), V(\gamma))  \approx (B_n, {\bf
R}^n)$  and $(G(\gamma_0), V(\gamma_0))  \approx (D_n, {\bf R}^n)$ for some
$\gamma_0 \subset \Gamma_W$ with $n=n(\gamma_0) \ge 4$;

\noin iii) $(W, V) \approx (A_{\Delta}, {\bf R}^{\Delta}_0)$ in the other
cases.

\smallskip

 From now on we consider Coxeter graphs as weighted graphs. As usual, the
weight of an edge $(s', s'')$ is the order of the product $s's''$. Since $W$ is
a locally finite group, the weights on $\Gamma_W$ take only finite values.
Recall that the Coxeter graphs of types $A_n$ and $D_n$ have only edges of
weight 3, while in any of the Coxeter graphs of type $B_n$ there are edges of
weight 4.  Thus, if the assumption of (i) holds, then $(G(\gamma'), V(\gamma'))
\approx (B_{n'}, {\bf R}^{n'}) $ for any $\gamma' \supset \gamma_0$ with
$n'=n(\gamma') > 8$, and thus $(W, V)$ is the inductive limit of some net of
pairs $(B_n, {\bf R}^n)$.

Next we show that for $\gamma \subset \gamma'$, where $\gamma \supset \gamma_0$
and $n(\gamma) > 8$, the embedding  $(G(\gamma), V(\gamma)) \subset
(G(\gamma'), V(\gamma'))$ is coordinatewise. This means that under the
isomorphisms $(G(\gamma), V(\gamma)) \approx (B_n, {\bf R}^n) $ and
$(G(\gamma'), V(\gamma')) \approx (B_{n'}, {\bf R}^{n'}) $ the pair $(B_n, {\bf
R}^n)$ is a coordinate subpair of $(B_{n'}, {\bf R}^{n'}) $. Indeed, identify
$(G(\gamma'), V(\gamma'))$ with $(B_{n'}, {\bf R}^{n'}) $. Then $V(\gamma)
\subset {\bf R}^{n'}$ is spanned by a subsystem of the root system $(\pm
\epsilon_i,\,\, \pm \epsilon_i \pm  \epsilon_j), \,\,1\le i<j\le n'$, and
$G(\gamma)$ is generated by the corresponding orthogonal reflexions. A plane in
${\bf R}^{n'} $ spanned by roots may contain reflexion root vectors of  2, 3 or
4  different reflexions from $B_{n'}$. It is a coordinate plane precisely when
it contains 4 reflexions. Since $(G(\gamma), V(\gamma)) \approx (B_n, {\bf
R}^n) $, the subspace $V(\gamma)$ contains $\left( \begin{array}{cc} n \\ 2
\end{array} \right)$ planes of the latter type which span it. At the same time,
these planes should be coordinate planes of ${\bf R}^{n'}$. Therefore,
$V(\gamma)$ is a coordinate subspace ${\bf R}^n \subset {\bf R}^{n'}$, and so
$G(\gamma)$ coincides with the group $B_n$ generated by reflexions along those
roots of the above root system which belong to $V(\gamma)$.

For any of the graphs $\gamma \supset \gamma_0$ with $n=n(\gamma) > 8$ all the
vertices of the full Coxeter graph $\bar\gamma$ (see the notation above) are
divided in two types: those which correspond to sign change reflexions, i.e.
reflexions along the roots of the form $\epsilon_i, i=1,\dots,n,$ and others.
Being coordinatewise, embeddings of pairs respect this division.  Thus, it is
well defined in the inductive limit $\Gamma_W$. Note that the vertices of
change sign type in $\Gamma_W$ are those which are incident only with edges of
weight 4. Denote by $\Delta$ the set of all the vertices of $\Gamma_W$ of
change sign type. Fix $\delta_0 \in \Delta \cap \gamma_0$, and let $\epsilon_0
= \epsilon_{\delta_0}$ be one of the two opposite roots in $V(\gamma_0)$ which
correspond to the reflexion $\delta_0$. It is easily seen that for any $\gamma
\supset \gamma_0$ with $n=n(\gamma) > 8$ the orbit $G(\gamma)\,\epsilon_0$
consists of the roots of coordinate type $\pm \epsilon_i$ in $V(\gamma)$, and
so the class of conjugates of $\delta_0$ in $W$ coinsides with  $\Delta$.
Choosing one of any two opposite root vectors in the orbit $W(\epsilon_0)$ we
obtain a Hamel basis of the  $W$-invariant subspace ${\rm span}(W(\epsilon_0))$
which coinsides with $V$ since $W$ is assumed to be irreducible in $V$.  Thus,
we obtain a Hamel basis of $V$ formed by roots of coordinate type. This yields
an isomorphism $V \approx {\bf R}^{\Delta}$. The root system of $W$, which
consists of the vectors of the two orbits $W(\epsilon_0)$ and $W\,(\epsilon_0 +
\epsilon_1)$, where $\epsilon_1 \not = \epsilon_0$ is another coordinate vector
in $V(\gamma_0)$, corresponds under this isomorphism to the root system
$(\pm\epsilon_{\delta},\,\,\pm\epsilon_{\delta}\pm\epsilon_{\delta'}\,\vert\,\delta,\delta' \in \Delta, \delta \not = \delta')$ of the group $B_{\Delta}$. Therefore, $(W, V) \approx (B_{\Delta}, {\bf R}^{\Delta})$. This proves (i).

Next we consider the case (ii), where there is no subgroup $G(\gamma)\subset W$
of type $B_n$, but at least one of them, say $G(\gamma_0)$, has type $D_n$ for
some $n=n(\gamma_0)\ge 4$.  First we show that any  subgroup $G(\gamma)\supset
G(\gamma_0)$ is of type $D_{n(\gamma)}$, and all the embeddings $G(\gamma)
\hookrightarrow G(\gamma')$ are coordinatewise.

The group $D_4$ contains 4 pairwise commuting reflexions along the root vectors
$\epsilon_1\pm\epsilon_2, \epsilon_3\pm\epsilon_4$. If $v_1,\dots, v_4$ are 4
mutually orthogonal root vectors from the root system $(\pm(\epsilon_i -
\epsilon_j)\,|\,1\le i<j\le n')$ of type $A_{n'}$ and $L = {\rm
span}\,(v_1,\dots, v_4)$, then the only reflexions in $A_{n'}\,\vert\,L$ are
the orthogonal reflexions along $v_1,\dots, v_4$, and so $A_{n'}\,\vert\,L$
does not contain $D_4$. Therefore, the Coxeter group $G(\gamma)\supset
G(\gamma_0)$ is not of type $A_{n(\gamma)}$, and thus it must be of type
$D_{n(\gamma)}$.

Let F be a subspace of dimension 4 of ${\bf R}^{n'}$ generated by 4 mutually
orthogonal root vectors from the root system $(\pm\epsilon_i \pm
\epsilon_j),\,1\le i<j\le n'$, of type $D_{n'}$, and $G(F)\subset D_{n'}$ be
the subgroup generated by the orthogonal reflexions along the roots in $F$.
Then $G(F)$ is irreducible (and of type $D_4$) iff $F$ is a coordinate subspace
of ${\bf R}^{n'}$. Thus, if $(G(\gamma), V(\gamma)) \subset (D_{n'}, {\bf
R}^{n'})$ is of type $D_n$, where $n=n(\gamma) \ge 4$, then $V(\gamma)$ is a
coordinate subspace of $ {\bf R}^{n'}$.

Fix a reflexion vector $v_0$ of a reflexion from $G(\gamma_0)$. If
$G(\gamma)\supset G(\gamma_0)$, then the orbit $G(\gamma)\,v_0$ is a root
system of type $D_{n(\gamma)}$ of the Coxeter group $G(\gamma)$. Consider the
infinite root system $W\,(v_0)$ in $V$. Since $W$ is irreducible, this system
is complete in $V$. Note that two roots $v', v''$ of $D_n$ are contained in the
same coordinate plane in ${\bf R}^n$ iff the sets of their neigborhooding
vertices in the full Coxeter graph $\Gamma_{D_n}$ coincide. In this case the
vectors $({{\pm v'\pm  v''}\over 2}) = (\pm\epsilon_i,\, \pm \epsilon_j),
\,i\not = j$,  are contained in the coordinate axes which are the intersections
of coordinate planes. The same pairing  is defined on the above root system of
$W$. In this way, fixing one of any two opposite vectors $({{\pm v'\pm
v''}\over 2})$ arbitrarily, we obtain a Hamel basis $\Delta$ in $V$, which in
turn provides us with an isomorphism $(W, V) \approx (D_{\Delta}, {\bf
R}^{\Delta})$. This proves (ii).

Assume further that any subgroup $G(\gamma')\subset W$ with $n'=n(\gamma') > 8$
is of type $A_{n'}$, where $A_{n'}$ acts by permutations in $ {\bf
R}^{n'+1}_0$. Let $(G(\gamma), V(\gamma)) \subset (A_{n'}, {\bf R}^{n'+1}_0)$
be of type $A_n$. We will show that $V(\gamma)$ is a coordinate subspace of
${\bf R}^{n'+1}_0$. Let $s_{ij}$ be the orthogonal reflexions (transpositions)
along the roots $\pm(\epsilon_i - \epsilon_j)\,\,,1 \le i <j \le n'+1$. Put $$I
= \{i \in \{1,\dots, n'+1\}\,\vert\,s_{ij} \in G(\gamma)\,\,\, {\rm for\,\,\,
some}\,\,\, j \in \{1,\dots, n'+1\}\}\,\,\,.$$  Thus, if $s_{kl} \in
G(\gamma)$, then $k,l \in I$. Vice versa, $s_{kl} \in  G(\gamma)$ for any pair
$k, l \in I, k \not = l$. This follows from the connectedness of the Coxeter
graph $\Gamma_{G(\gamma)}$ and the following remark: if $s_{ij} \in  G(\gamma)$
and $s_{jk} \in  G(\gamma)$, then $s_{ik} \in  G(\gamma)$. Indeed,
$s_{jk}(\epsilon_i - \epsilon_j) = \epsilon_i - \epsilon_k$ and thus
$s_{ik}=s_{jk} s_{ij} s_{jk}$. Now we see that $V(\gamma) = {\bf R}^I_0 = {\rm
span}\,(\epsilon_i - \epsilon_j \,\vert\,i,j \in I)$ is a coordinate subspace
of ${\bf R}^{n'+1}_0$, and $G(\gamma)\subset A_{n'}$ is a subgroup of
permutations of the set $I$.

On the set of edges of the full Coxeter graph  $\Gamma_{A_n}$ consider the
following equivalence relation: $(s_{i_1,j_1}, s_{i_2,j_2}) \sim (s_{k_1,l_1},
s_{k_2,l_2})$ iff these four transpositions have an index in  common. Then this
index is the same for the whole equivalence class, so that the set of classes
is $\{1,\dots,n\}$. Since this equivalence relation is compatible with the
embedding of pairs $(G(\gamma), V(\gamma)) \hookrightarrow (G(\gamma'),
V(\gamma'))$, it can be defined as well on the whole graph $\Gamma_W$. Let
$\Delta$ be the set of the equivalence classes. Let $v \in \Gamma_W$ be a
vertex. Then all the edges incident with $v$ belong to two  different classes
$\delta, \delta' \in \Delta$, where each class $\delta \in \Delta$ consists of
the edges of a complete subgraph of $\Gamma_W$, and each pair of these complete
subgraphs which correspond to some $\delta, \delta' \in \Delta, \,\,\delta \not
= \delta',$ has exactly one vertex $v(\delta, \delta')$ in common. It is easily
seen that the action of $W$ on $\Gamma_W$ by inner automorphisms is locally
finite and compatible with the equivalence relation, and so it induces the
action of $W$ on $\Delta$ by finite permutations, such that the reflexions in
$W$ act as transpositions.

Fix a reflexion $s_0 = s(\delta, \delta') \in W$ which corresponds to the
transposition $(\delta, \delta')$, with the reflexion vector $e_0 = e(\delta,
\delta')$. Then the orbit $W\,(e_0)$ is a root system of $W$ which spans $V$.
Fixing one of the each two opposite roots, we obtain a Hamel basis of $V$ which
corresponds to the basis of ${\bf R}^{\Delta}_0$ consisting of the root vectors
$\epsilon_{\delta} - \epsilon_{\delta'} \,(\delta, \delta' \in \Delta, \delta
\not = \delta')$ of $A_{\Delta}$. This gives an isomorphism $(W, V) \approx
(A_{\Delta}, {\bf R}^{\Delta}_0)$. The proof is  complete. \,\,\,$\bigcirc$

\section{Total families of isometric reflexions}

Denote by ${\rm {\rm IR}}(E)$ the set of all the reflexions in ${\rm Iso}\, E$,
and let $W = W (E) $ be the subgroup of ${\rm Iso}\, E$ generated by the
reflexions from ${\rm {\rm IR}}(E)$. In this section we assume that ${\rm {\rm
IR}}(E)$ contains a total subset of reflexions $\{ s_{\alpha} = s_{e_{\alpha},
e^*_{\alpha}}\}_{\alpha \in A}$, which means that the family of linear
functionals $T = \{e^*_{\alpha}\}_{\alpha \in A} \subset E^*$ is a total
family.

\bigskip

\noin {\bf 4.1. Lemma.} {\sl Let $g_1 , g_2 \in {\rm Iso}\, E$. In the
notation as above assume that  $g_1(e_{\alpha}) = g_2(e_{\alpha})$ for all
$\alpha \in A$. Then $g_1 = g_2$.} \bigskip

\noin {\sl Proof}. Put $g_0 = g_1^{-1} g_2$. Then $g_0(e_{\alpha}) =
e_{\alpha}$ for all ${\alpha} \in A$. Since $s_{\alpha}$ is the only isometric
reflexion in the direction of $e_{\alpha}$ (see Remark 2.2.b), it coincides
with $s'_{\alpha} = g_0 s_{\alpha} g_0^{-1} = s_{e_{\alpha}, g_0^{*
-1}(e^*_{\alpha})}$, and so $ g_0^{* -1}(e^*_{\alpha}) =  e^*_{\alpha}$, i.e. $
g_0^* (e^*_{\alpha}) =  e^*_{\alpha}$ or, in other words, $e^*_{\alpha}(g_0(v)
- v) = 0$ for all $\alpha \in A$. Since $T$ is total, it follows that $g_0 =
1_E$. \,\,\,$\bigcirc$\\

 Let, as before, $G_0$ be the strong identity component of ${\rm Iso}\, E$, and
let $W$ be the group generated by the reflexions in ${\rm Iso}\, E$. \bigskip

\noin {\bf 4.2. Lemma.} {\sl  $W$ is locally finite iff $G_0$ is trivial.
}\bigskip

\noin {\sl Proof}. Assume that $G_0$ is trivial. To prove that $W$ is locally
finite it is enough to show that each subgroup $W'$ of $W$ generated by a
finite number of reflexions $\{s_i = s_{e_i, e^*_i}\}_{i=1,\dots,n} \subset
{\rm {\rm IR}}(E)$ is finite. Suppose that $W'$ is an infinite group. Put $F' =
{\rm span}(e_i\,|\, i=1,\dots,n)$.  Let $G'$ be the closure of $W'$ in ${\rm
Iso}\, E$ in the strong operator topology. It is easily seen that the closed
subspace $M' = \bigcap _{i=1}^n {\rm Ker}\,e^*_i$ is a complementary subspace
of $F'$, i.e. $E = M' \oplus F'$, and it coincides with the fixed point
subspace of $W'$. Hence, it also coincides with the fixed point subspace of
$G'$. It follows that $G' = 1_{M'} \oplus {\bar G}' $, where ${\bar G}' \subset
O(F')$ is the closure of $W' \, \vert \, F'$ in
${\rm Iso}\, F'$. Thus, $G'$ is a compact Lie group, and being infinite it has
a non-trivial identity component. This is a contradiction.

Assume now that $G_0$ is non-trivial. Then, as it was shown in the proof of
Proposition 2.6, there exist reflexions $s' , s'' \in {\rm {\rm IR}}(E)$ such
that the angle $\alpha (s' , s'')$ is irrational modulo $\pi$, and so the
subgroup of $W$ generated by these two reflexions is infinite.
\,\,\,$\bigcirc$\\

Remind that a group $G$ of operators in a Banach space $E$ is called {\it
topologically irreducible} if it has no nontrivial closed invariant subspace.
\bigskip

\noin {\bf 4.3. Lemma.} {\sl Let $W'$ be a group generated by a set of
reflexions $\{ s_{\alpha} = s_{e_{\alpha}, e^*_{\alpha}}\}_{\alpha \in A'}
\subset {\rm {\rm IR}}(E)$. Then $W'$ is topologically irreducible iff the
following two conditions are fulfilled:

\noin i) The system of vectors $(e_{\alpha}\, \vert\, \alpha \in A')$ is
complete, i.e. $E = {\overline {\rm span}}(e_{\alpha}\, \vert \alpha \in A')$.

\noin ii) The Coxeter graph $\Gamma_{W', A'}$ is connected. }

\bigskip

\noin {\sl Proof}. Since the closed subspace $E' = {\overline {\rm
span}}(e_{\alpha}\, \vert \,\alpha \in A')$ is invariant with respect to $W'$,
the first condition (i) is necessary for $W'$ being irreducible. Let $\Gamma '$
be a connected component of $\Gamma_{W', A'}$. It is easily seen that the
closed subspace $F' = {\overline {\rm span}}(e\, \vert \, s_{e, e^*} \in
\Gamma')$ is invariant, too. Thus, the second condition (ii) is necessary, too.

Assume further that (i) and (ii) are fulfilled. Let $F'$ be a closed invariant
subspace of $W'$. Put $A = \{\alpha \in A'\, \vert \,e_{\alpha} \in F'\}$ and
$B = \{\alpha\in A'\, \vert \,e_{\alpha} \notin F'\}$. Being invariant $F'$ is
contained in ${\rm Ker}\,e^*_{\beta}$ for each  $\beta \in B$. It follows that
$\alpha (s_{\alpha} , s_{\beta}) = \pi /2$, and so $[s_{\alpha} , s_{\beta}] =
1_E$ for any $ \alpha\in A , \beta \in B$. By (ii) this implies that either $A
= \emptyset$ or $B = \emptyset$. From (i) it easily follows that the system of
linear functionals $(e^*_{\alpha}\,\vert\,\alpha\in A)\subset E^*$ is total.
Thus, if $A = \emptyset$, then $F' \subset \bigcap \{{\rm
Ker}\,e^*_{\alpha}\,\vert\,\alpha \in A'\} = \{\bar 0 \}$, and if $B =
\emptyset$, then $F' \supset {\overline {\rm span}}\,(e_{\alpha}\,\vert\,\alpha
\in A') = E$. In any case, $F'$ is not a proper subspace. This shows that $W'$
is topologically irreducible.\,\,\,$\bigcirc$ \\

Let things be as in Lemma 4.3. Consider the algebraic linear subspace $V' =
{\rm span}\,(e_{\alpha}\, \vert\, \alpha \in A')$. The group $W'$ is
algebraically irreducible in $V'$ iff the Coxeter graph $\Gamma_{W', A'}$ is
connected. In this case $W'$ is topologically irreducible in the closed
subspace $E' = {\overline V'} =  {\overline {\rm span}}\,(e_{\alpha}\, \vert
\,\alpha \in A')$. If $W'$ is finite and the Coxeter graph  $\Gamma_{W'}$ is
connected, then ${\rm \dim}\,V' = n < \infty$ and $W'\,\vert\,V'$ is a finite
Coxeter group, i.e. a finite irreducible group generated by orthogonal
reflexions in ${\bf R}^n$ (here we identify $V'$ with ${\bf R}^n$ by choosing
an orthonormal basis with respect to an invariant scalar product in $V'$). Let
the group ${\rm Iso} E$ be discrete in the strong operator topology, i.e. $G_0
= \{1_E\}$. Then by Lemma 4.3, $W'$ is a locally finite group. If ${\rm
\dim}\,V' = \infty$ and the Coxeter graph  $\Gamma_{W'}$ is connected, then
$W'$ is an infinite Coxeter group, and by Proposition 3.2, it is isomorphic to
one of the groups $A_{\Delta}, \,\,\,B_{\Delta}, \,\,\,D_{\Delta}$. The next
proposition shows that if the pair $(W', V')$ is maximal, it can not be of type
$D_{\Delta}$.

\bigskip

\noin {\bf 4.4. Proposition.} {\sl Let the notation be as above. If ${\rm
\dim}\,V' = \infty$ and $(W', V') \approx (D_{\Delta}, {\bf R}^{\Delta})$, then
the group $W'$ can be extended to a subgroup $W'' \subset {\rm Iso} E$
generated by reflexions along vectors in $V'$ and such that $(W'', V') \approx
(B_{\Delta}, {\bf R}^{\Delta})$}.

\bigskip

For the proof we need the following lemma on partial orthogonal decompositions
in Banach spaces.

\bigskip

\noin {\bf 4.5. Lemma.} {\sl Let $\{p_i\}_{i=1,2,\dots}$ be a sequence of
projections in a Banach space $E$ such that

\smallskip

\noin a)\,\,\, $1_E - 2p_i \in {\rm Iso} E \,\,\,{\rm for\,\,\, all}\,\,\,
i=1,2,\dots$;

\noin b) \,\,\, the projections $p_i$ are mutually orthogonal, i.e. $p_i p_j =
0\,\,\, {\rm for\,\,\, all}\,\,\, i \ne j$.

\smallskip

\noin Then $\,\,\,\limsup_{i \rightarrow \infty} ||(1_E - p_i)(x)||_E = ||x||_E
\,\,\,{\rm for\,\,\, all}\,\,\,x \in E$.}

\bigskip

\noin {\sl Proof}. By (a), we have $||p_i||_E = ||1_E - p_i||_E = 1\,\,\,{\rm
for\,\,\, all}\,\,\,i=1,2,\dots$. From (a) and (b) it follows that
$\prod_{i=1}^k (1_E - 2p_i) = 1_E - 2\sum_{i=1}^k p_i \in {\rm Iso} E$, and so
$||\sum_{i=1}^k p_i||_E = ||1_E - \sum_{i=1}^k p_i||_E = 1$, as well.

Assume that there exist $x_0 \in E$ and $\epsilon_0 > 0$ such that $$||(1_E -
p_i) (x_0)||_E \le ||x_0||_E - \epsilon_0 \,\,\,{\rm for\,\,\, all}\,\,\,
i=1,2,\dots\,\,\,.$$ Then $$ \Vert {1\over k}\sum_{i=1}^k (1_E - p_i) (x_0)
\Vert_E \le {1\over k} \sum_{i=1}^k \Vert (1_E - p_i)(x_0)\Vert_E \le ||x_0||_E
- \epsilon_0 \,\,\,.$$
Therefore
$$\epsilon_0\le ||x_0||_E - \Vert {1\over k}\sum_{i=1}^k (1_E - p_i) (x_0)
\Vert_E \le \Vert x_0 - {1\over k}\sum_{i=1}^k (1_E - p_i) (x_0)\Vert_E =
{1\over k}\Vert (\sum_{i=1}^k  p_i) (x_0)\Vert_E \le  {1\over k}||x_0||_E\,\,\,
.$$
This is a contradiction. $\,\,\,\bigcirc$\\

\bigskip

\noin {\sl Proof of Proposition 4.4}. Identify $V'$ with ${\bf R}^{\Delta}$
via an isomorphism $(W', V') \approx (D_{\Delta}, {\bf R}^{\Delta})$, and
consider in $V'$ the root system $\{\pm \epsilon_{\delta'} \pm
\epsilon_{\delta''}\}$
of type $D_{\Delta}$.  Denote by $s_{\delta', \delta''}^+$ the isometric
reflexion along the vector $v_{\delta' ,\delta''}^+ = \epsilon_{\delta'} +
\epsilon_{\delta''}, \delta', \delta'' \in \Delta, \delta' \ne \delta'',$  and
by $s_{\delta', \delta''}^-$ the isometric reflexion along the vector
$v_{\delta', \delta''}^- = \epsilon_{\delta'} -  \epsilon_{\delta''}$. Put $d
_{\delta', \delta''} = s_{\delta', \delta''}^+s_{\delta', \delta''}^- \in W'$,
so that $d _{\delta', \delta''}$ is the operator of change of signs of the
coordinates $\delta'$ and $\delta''$.

Choose a countable subset $\{\delta_i\}_{i=1,2,\dots} \subset \Delta$ and put
$d_{i,j} = d _{\delta_i, \delta_j}$. Then the involutions $d_{i,j}$ pairwise
commute and $d_{n,k}d_{k,m}=d_{n,m}$. The orthogonal projections $p_{i,j} =
{1\over 2}(1_E - d_{i,j})$ onto planes also pairwise commute.
For each triple of different indices $n, m, k$ consider the one-dimensional
projection $p_{n}^{k,m} = p_{n,k}p_{n,m}$. Since  $p_{n}^{k,m}$ and
$p_{n}^{i,j}$ commute and have the same image, they coincide; indeed,
$p_{n}^{k,m} = p_{n}^{i,j}p_{n}^{k,m}=p_{n}^{k,m}p_{n}^{i,j}=p_{n}^{i,j}$.
Denote by $p_n$ their common value, and consider the corresponding reflexion
$s_n = 1_E - 2p_n$ along the coordinate vector $\epsilon_{\delta_n}$. It is
easily seen that $s_{n} s_{m} = d_{n,m}$ and
$s_m (1_E - p_{m, k}) = 1_E - p_{m, k}$. By Lemma 4.5, for a
fixed $n \in {\bf N}$
and for any $x \in E, \epsilon > 0$ there exist $k, m  \in {\bf N}$ such that
$$\Vert s_n (x)\Vert_E \le \Vert (1_E - p_{k, m})s_n  (x)\Vert_E + \epsilon =
\Vert s_n (1_E - p_{k, m})(x) \Vert_E + \epsilon = $$
$$\Vert s_n s_m (1_E - p_{k, m}) (x) \Vert_E + \epsilon = \Vert d_{n,m} (1_E -
p_{k, m}) (x) \Vert_E + \epsilon = $$
$$\Vert (1_E - p_{k, m}) (x) \Vert_E + \epsilon \le \Vert x \Vert_E +
\epsilon\,\,\,.$$ It follows that $s_n \in {\rm Iso} E$. Since $\delta_n \in
\Delta$ is taken as arbitrary, this implies that for any $\delta \in \Delta$
there exists an isometric reflexion along the vector $\epsilon_{\delta}$.
Thus, the group ${\rm Iso} E$ contains the subgroup $W''$ generated by
reflexions along vectors of the root system $\{\pm \epsilon_{\delta}, \pm
\epsilon_{\delta'} \pm \epsilon_{\delta''}\}$ of type $B_{\Delta}$.
$\,\,\,\bigcirc$\\

\bigskip

 \noin {\bf 4.6. Corollary.} {\sl Let $\dim E = \infty$, $G_0 = G_0 (E) =
\{1_E\}$, and the group $W=W(E)$ generated by all the isometric reflexions in
$E$ be topologically irreducible. Then $W$ is an infinite Coxeter group of type
$A_{\Delta}$ or $B_{\Delta}$.} \\

\bigskip

 \noin {\bf 4.7.} {\sl Remark}. Let $\dim E = n < \infty$. Then Proposition 4.4
still holds in the case when $n$ is odd. Indeed, in this case $B_n = W''$ is
the subgroup of the group ${\rm Iso} E$ generated by the subgroup $W' = D_n$
and the element $-1_E$. But for $n$ even the statement of Proposition 4.4 in
general is not valid. As an example, consider $E={\bf R}^n$, where $n = 2k \ge
4$, with the unit ball $B(E)$ being the convex hull of the $D_n$-orbit of the
point $v_0 = (1,2,\dots,n)$. Then the image $s_n (v_0) = (1,2,\dots,n-1,-n)$ of
$v_0$ by the reflexion $s_n = s_{\epsilon_n, \epsilon_{n}^*}$ does not belong
to $B(E)$ (indeed, it is separated from $B(E)$ by the hyperplane $-x_n +
\sum_{i=1}^{n-1} x_i = {n(n+1)\over 2}$). Hence $B(E)$ is not invariant with
respect to the action of the Coxeter group $B_n$ on ${\bf R}^n = E$, and so
$B_n$ is not a subgroup of ${\rm Iso}\, E$.

\section{ Hilbert and Coxeter decompositions}

Let, as before, $E$ be a Banach space with a total family of isometric
reflexions. In this section we construct a partial orthogonal decomposition of
$E$ which consists of two parts: {\it Hilbert decomposition} into a direct sum
of biorthogonally complemented Hilbert subspaces, and {it Coxeter
decomposition} into a direct sum of closed subspaces endowed with topologically
irreducible Coxeter groups generated by isometric reflexions. In a sense, this
decomposition is orthogonal (see Lemma 5.4 and Proposition 5.6). Both of these
decompositions are stable under the action of the isometry group ${\rm Iso} E$,
and the second one is fixed under the action of the identity component $G_0$.
The main result of the section, Theorem 5.7, is a kind of a structure theorem
for the isometry group ${\rm Iso} E$.\\

\bigskip

 \noin {\bf 5.1.} {\it  Notation.} As above, by ${\rm IR}(E)$ we denote  the
set of all the isometric reflexions in $E$ which is assumed to be total. To
each subspace $V$ of $ E$ we attach two closed subspaces, {\sl  the  kernel}
$$V_0 = {\overline {\rm span}} \,(e \in V \,|\,s_{e, e^*} \in {\rm IR} (E)\,)$$
and {\sl the  hull} $$\hat {V} = \bigcap \{ {\rm Ker}\, e^*\,|\,s_{e, e^*} \in
{\rm IR} (E)\,, \,V \subset {\rm Ker}\, e^*\}\,\,;$$ we put $\hat {V} = E$ if
there is no $s_{e, e^*} \in {\rm IR} (E)$ such that $ V \subset {\rm Ker}\,
e^*$. It is easily seen that

\smallskip

\noin {\it i) $V_0 \subset {\overline V} \subset \hat V$,

\noin ii) $V_{00} = V_0\,,\,\hat {\hat {V}} = \hat {V}$, and

\noin iii) if $V \subset V'$, then $V_0 \subset V'_0$ and $\hat V \subset \hat
{V'}$.}

\smallskip

\noin Observe that possibly $V_0$ resp. ${\overline V}$ is a proper subspace of
${\overline V}$ resp. $\hat {V}$. For instance, this is the case when $E = {\rm
l}_{\infty}$ and $V = {\rm c}$ (the subspace of convergent sequences); indeed,
then $\hat {V} = E$ and $V_0 = {\rm c}_0$.

Denote also ${\rm IR}_V = \{s_{e, e^*} \in {\rm IR} (E)\,|\,e \in V\}$. Let
$W_V$ be the group generated by the reflexions from ${\rm IR}_V$. \bigskip

\noin {\bf 5.2. Coxeter decomposition.} This is a partial subspace
decomposition defined on the fixed point subspace $F = {\rm Fix} \, G_0$ of the
group $G_0 = G_0(E)$. Let $\Gamma _F = \Gamma_{W_F}$ be the full Coxeter graph
of the group $W_F$, and let $  \cal A$ be the set of the connected components
of $\Gamma_F$. For $\alpha \in \cal A$ denote by ${\rm IR}_{\alpha}$ the set of
reflexions in ${\rm IR}_F$ which correspond to vertices of the component
$\alpha$ of $\Gamma_F$. Put $V_{\alpha} = {\overline {\rm span}}\,(e\,|\,s_{e,
e^*} \in {\rm IR}_{\alpha})$; so, ${\rm IR}_{\alpha} = {\rm IR}_{V_{\alpha}}$.
Put also $W_{\alpha} = W_{V_{\alpha}}$. Then $V_{\alpha}$ is a closed subspace
of the kernel $F_0$, and the group $W_{\alpha}\,|\,V_{\alpha}$ is topologically
irreducible. By the discussion after Lemma 4.3, $W_{\alpha}$ is a Coxeter
group. If $\dim V_{\alpha} = \infty$, then by Corollary 4.6, $W_{\alpha}$ has
type $A_{\Delta}$ or $ B_{\Delta}$.

The set $\cal A$ can be devided into equivalence classes which correspond to
the isomorphism types of the Coxeter pairs $(W_{\alpha}, V_{\alpha})$. Since
$G_0$ is a normal subgroup of the group ${\rm Iso}\,E$, its fixed point
subspace $F$ is invariant with respect to ${\rm Iso}\,E$; the same is true for
the kernel $F_0$ and the hull $\hat {F}$. Each isometry $g \in {\rm Iso}\,E$
acts (by conjugation) on the set ${\rm IR}_F$ and also on the graph $\Gamma_F$,
and so on the set $\cal A$. It is clear that the above partition of $\cal A$ is
stable under this action and its equivalence classes are invariant.\bigskip

\noin {\bf 5.3. Hilbert decomposition.} Consider the following equivalence
relation defined on the set ${\rm IR}(E) \setminus {\rm IR}_F\,$: $$s_{e, e^*}
\sim s_{e', e'^*}\,\,\, {\rm iff} \,\,\, e' \in G_0\, e\,\,\,.$$ Let $\cal B$
be the set of its equivalence classes. By Theorem 1, to each $\beta \in \cal
{B}$ there corresponds the unique Hilbert subspace $H_{\beta} = {\overline {\rm
span}}\,(G_0 \,e\,|\,s_{e, e^*} \in \beta)$ and the unique bicontractive
projection $p_{\beta}: E \rightarrow H_{\beta}$ satisfying all the properties
of Theorem 1.b. An isometry $g \in {\rm Iso}\,E$ induces the action $g_*$ on
the set $\cal B$ which is defined as follows: $g_* \beta = \beta '$ iff
$g(H_{\beta}) = H_{\beta'}$. In particular, the orthogonal bases in $H_{\beta}$
and in $H_{\beta'}$ are of the same cardinality. The following lemma shows that
this partial decomposition into Hilbert subspaces is orthogonal; moreover, all
of the subspaces $H_{\beta}$ are orthogonal to the fixed point subspace $F$.

Let ${\rm IR}_{\beta} = {\rm IR}_{H_{\beta}}$. Note that ${\rm IR} (E) = {\rm
IR}_{\cal {A}} \cup {\rm IR}_{\cal {B}}$, where ${\rm IR}_{\cal {A}} = {\rm
IR}_F = \bigcup_{\alpha \in {\cal {A}}} {\rm IR}_{\alpha}$ and ${\rm IR}_{\cal
{B}} =  \bigcup_{\beta \in \cal {B}} {\rm IR}_{\beta}$. \\

\bigskip

\noin {\bf 5.4. Lemma.} {\sl a) Let $s, s' \in {\rm IR}(E)$. If $[s, s'] \neq
1_E$, then $s, s'$ belong either to the same subset ${\rm IR}_{\alpha}$, where
$\alpha \in \cal {A}$, or to the same subset ${\rm IR}_{\beta}$, where $\beta
\in \cal {B}$.

\noin b)  The projection $p_{\beta}$ commutes with any reflexion $s \in {\rm
IR}(E)$ for any  $\beta \in \cal {B}$.

\noin c) Furthermore, $ p_{\beta}p_{\beta'} = 0$ for any $\beta, \beta' \in
\cal {B}, \beta \neq \beta'$, and $p_{\beta}\,|\,F = 0$ for any  $\beta \in
\cal {B}$.}

\bigskip

\noin {\sl Proof. a}. Let $s_i = s_{e_i, e_i ^*}, \,i=1,2$, be two arbitrary
distinct reflections from $ {\rm IR}_{\beta}$, where $\beta \in \cal {B}$.
Being restricted to the Hilbert subspace $H_{\beta}$ the rotation $r = s_1 s_2
\in {\rm Iso}\,E$ in the plane $L = {\rm span}\,(e_1, e_2)$ belongs to the
connected component $G_0 (H_{\beta})$ of the orthogonal group, and so by
Theorem 1.c, $r \in G_0$. Since $F = {\rm Fix}\, G_0 \subset {\rm Fix}\,r =
{\rm Ker}\,e_1^* \cap {\rm Ker}\,e_2^* $, we have that $e_i^* (e) = 0 $ for
each $e \in F$. Therefore, if $s= s_{e, e^*} \in {\rm IR}_F$ then by Lemma 2.3,
$\alpha (s_i, s) = {\pi\over 2}$, and thus $[s_i, s] = 1_E,\,\, i=1,2$ (see
Remark 2.2.c).

If $\beta' \in \cal {B}$ and $\beta' \neq \beta$, then the subspace
$H_{\beta'}$ is invariant with respect to the rotation $r = s_1 s_2 \in G_0$.
One may assume that $r\,|\,L \neq -1_L$, and so either $H_{\beta'} \subset {\rm
Fix}\,r$ or $L \subset H_{\beta'}$. The second case is impossible (indeed,
otherwise by the construction, we would have $H_{\beta} = H_{\beta'}$, and so
$\beta = \beta'$). Thus, $H_{\beta'} \subset {\rm Fix}\,r = {\rm Ker}\,e_1^*
\cap {\rm Ker}\,e_2^*$. As above, it follows that $[s_i, s] = 1_E$ for each $s
\in {\rm IR}_{\beta}$.

To prove (a) it remains to note that the definition of the set $\cal A$ (5.2)
yields that $[s, s'] = 1_E$ if $s \in {\rm IR}_{\alpha},\, s' \in {\rm
IR}_{\alpha'}\,$, where $\alpha, \alpha' \in \cal {A}$ and $\alpha \neq
\alpha'$.

Now, (b) and (c) easily follow from (a) by the construction of the projections
$p_{\beta}$ as in 2.8. $\,\,\,\bigcirc$

\bigskip

 \noin {\bf 5.5. Lemma.} {\sl a) $\hat {F}_0 = \hat {F} = F$.

\noin b) $(H_{\beta})_0 = \hat {H_{\beta}} = H_{\beta} $ for any $\beta \in
\cal {B}$.}

\bigskip

\noin {\sl Proof. a}. Let $s_{e, e^*} \in {\rm IR} (E)$ be such that $F_0
\subset {\rm Ker}\,e^*$. Then $e \notin F_0$, and hence $e \in H_{\beta}$ for
some $\beta \in \cal {B}$. As in the proof of Lemma 5.4, it follows that $F
\subset {\rm Ker}\,e^*$, and so $F \subset \hat {F_0} = \bigcap_{\beta \in \cal
{B}} {\rm Ker}\,p_{\beta}$. Since the subspace $\bigcup_{\beta \in \cal {B}}
H_{\beta} $ is $G_0$-invariant, it is clear that $\hat {F_0}$ is
$G_0$-invariant, too.

If $\hat {F_0} \neq F$, then there exists $g_0 \in G_0$ such that $g_0\,|\,\hat
{F_0} \neq 1_{F_0}$, and so $g_0(x) \neq x$ for some $x \in \hat {F_0}$. Note
that both $g_0 (x)$ and $x$ belong to ${\rm Ker}\,e^* $ for each $e^*$ such
that $s_{e, e^*} \in {\rm IR} (E) \setminus {\rm IR}_F$ (indeed, in this case
$g_0 (e) = e$, and thus $g_0^* (e^*) = e^*$). Therefore, $e^* (g_0 (x) - x) =
0$ for each $e^*$ as above, and also for each $e^*$ such that $s_{e, e^*} \in
{\rm IR} (E)$. Since the system of functionals $(e^*\,|\,s_{e, e^*} \in {\rm
IR} (E))$ is total, we have $g_0 (x) - x = 0$, which is a contradiction. This
proves (a).

\noin b. If $\hat {H_{\beta}} \neq H_{\beta}$ for some $\beta \in \cal {B}$,
then $(1_E - p_{\beta}) (x) \neq 0$ for some vector $x \in \hat {H_{\beta}}$.
By Lemma 5.4.b, the projection $p_{\beta}$ commutes with any reflexion $s  =
s_{e, e^*}  \in {\rm IR}(E)$. Thus, if $H_{\beta} \subset {\rm Ker}\,e^* $,
then also $\hat {H_{\beta}} \subset {\rm Ker}\,e^*$.  Therefore, $s(y) = y$ for
all $y \in \hat {H_{\beta}}\,,\, p_{\beta} s(y) = s p_{\beta}(y) =
p_{\beta}(y)$ and $s(1_E - p_{\beta})(y) = (1_E - p_{\beta})(y)$. The latter
means  that $(1_E - p_{\beta})(y) \in {\rm Ker}\,e^*$. Hence, $(1_E -
p_{\beta})(\hat {H_{\beta}}) \subset \hat {H_{\beta}}$.

 Now, we have $(1_E - p_{\beta})(x) \in{\rm Ker}\,e^*$ for any $e^*$ such that
$s_{e, e^*} \in {\rm IR} (E)$. This contradicts to the assumption that the
system ${\rm IR}(E)$ is total, since $(1_E - p_{\beta})(x) \neq 0$.
$\,\,\,\bigcirc$

\bigskip

 Put $R_0 = {\overline {\rm span}}\, (\bigcup_{\beta \in \cal B} H_{\beta})$
and $\hat R = \hat R_0$.

\bigskip

\noin {\bf 5.6. Proposition.}
{\sl  a) The subspace  $R_0 \dot{+}  F$ is closed, and if $p_{R_0 , F} : R_0 \dot{+}
F \to R_0$  is the first projection, then $1_{R_0 \dot{+}   F} - 2p_{R_0 , F} \in
{\rm Iso}\,(R_0 \dot{+}   F)$. Therefore, the projection $p_{R_0 , F}$ is
bicontractive.

\smallskip

\noin b) For any $\alpha \in {\cal A}$ there exists a projection $p_{\alpha}:
F_0 \dot{+}   {\hat R} \to V_{\alpha}$ such that

\noin i) $p_{\alpha}$ commutes with any reflexion $s \in {\rm IR}\,(E)$ and
$p_{\alpha}p_{\alpha'}=p_{\alpha'}p_{\alpha}=0$ resp.
$p_{\alpha}p_{\beta}=p_{\beta}p_{\alpha}=0$ for all $\alpha' \in {\cal
A}\,,\,\alpha' \neq \alpha\,,\,\beta \in {\cal B}$;

\noin ii) $||p_{\alpha}||_{F_0 \dot{+}   {\hat R}} \le 2$ and $||1_{F_0 \dot{+}   {\hat
R}} - p_{\alpha}||_{F_0 \dot{+}   {\hat R}} = 1$, if the latter projection is
non-zero;

\noin iii) moreover, if the Coxeter group $W_{\alpha}$ is a group of type
$B_{\Delta}$, then $1_{F_0 \dot{+}   \hat R} - 2p_{\alpha} \in {\rm Iso}\,(F_0 \dot{+}
 {\hat R})$, and so $||p_{\alpha}||_{F_0 \dot{+}   {\hat R}} = 1$, too.

\smallskip

\noin c) The subspace $F_0 \dot{+}   {\hat R}$ is closed, and if both subspaces
$F_0$ and $\hat R$ are non-trivial and $p_{F_0 , {\hat R}} : F_0 \dot{+}   {\hat R}
\to F_0$  is the first projection, then $||1_{F_0 \dot{+}   {\hat R}} - p_{F_0 ,
{\hat R}}||_{F_0 \dot{+}   {\hat R}} = 1$ and $||p_{F_0 , {\hat R}}||_{F_0 \dot{+}
{\hat R}} \le 2$.}

\bigskip

\noin {\sl Proof. a}. Let $x = x_1 + x_2$, where $x_1 \in R_0$ and $x_2 \in F$.
For any $\epsilon > 0$ there exists a finite subset $\sigma \subset {\cal B}
$ and a vector $x_1^{\sigma} \in \bigoplus\limits_{\beta \in \sigma} H_{\beta}$
such that $||x_1 - x_1^{\sigma}||_E < \epsilon$. Since $u_{\sigma} =
\prod_{\beta \in \sigma} (1_E - 2p_{\beta}) \in {\rm Iso}\,E$ and
$u_{\sigma}(x_1^{\sigma}) = -x_1^{\sigma}\,,\,u_{\sigma}(x_2) = x_2$, we have
$||x_1^{\sigma} + x_2||_E = ||- x_1^{\sigma} + x_2||_E$. Thus, if $R_0 \neq
\{0\}$ and $F \neq \{0\}$, then $||1_{R_0 \dot{+}   F} - 2p_{R_0 , F}||_{R_0 \dot{+}
F} = 1$, and therefore also $||p_{R_0 , F}||_{R_0 \dot{+}   F} = ||1_{R_0 \dot{+}   F}
- p_{R_0 , F}||_{R_0 \dot{+}   F} = 1$, if both subspaces $R_0$ and $F$ are
non-trivial. By the closed graph theorem, this implies that the subspace $R_0
\dot{+}   F$ is closed.

\noin $b$. If $\dim V_{\alpha} < \infty$, put $p'_{\alpha} = (1/{\rm
card}\,W_{\alpha})\sum _{g \in W_{\alpha}} g$. Then $p'_{\alpha}$ is a
projection on the fixed point subspace $F_{\alpha}$ of the group $W_{\alpha}$,
which coincides with $\cap_{s_{e, e^*} \in {\rm IR}_{\alpha}} {\rm Ker}\,e^*$,
and thus it is a complementary subspace to $V_{\alpha} = {\rm Ker}\,
p'_{\alpha}$. It is clear that $||p'_{\alpha}||_E = 1$, and so $||1_E -
p'_{\alpha}||_E \le 2$. From Lemma 5.4 and the definition of $p'_{\alpha}$ it
follows that the projection $p_{\alpha} = (1_E - p'_{\alpha})\,|\,(R_0 \dot{+}
F)$ satisfies (i); by the above inequalities, it also satisfies  (ii).

Next consider the case when $\dim V_{\alpha} = \infty$ and the Coxeter group
$W_{\alpha}$ is of type $A_{\Delta}$. For a finite subset $\sigma \subset
\Delta$ denote by $V_{\sigma}$ the subspace generated by the root vectors
$\epsilon_{\delta} - \epsilon_{\delta'}$, where $\delta\,,\,\delta' \in
\sigma\,,\,\delta\neq\delta'$, and by $W_{\sigma}$ the Coxeter group of type
$A_n$, where $n = \dim V_{\sigma}$, generated by the isometric reflexions along
these vectors. Define the projections $p'_{\sigma}$ resp. $p_{\sigma}$ in the
same way as $p'_{\alpha}$ resp. $p_{\alpha}$ above. It is clear that
$p_{\sigma}$ commutes with any reflexion $s \in {\rm IR}\,(E) \setminus {\rm
IR}_{\alpha}$ and satisfies all the other properties in (i), (ii).  It is
easily seen that the net $(p_{\sigma})$ is strongly convergent to the identity
on the subspace $ V_{\alpha}$, and that all the projections $p_{\sigma}$ vanish
on the subspace ${\hat R} \dot{+}   V'_{\alpha}$, where $V'_{\alpha} = {\overline
{\rm span}}\,(\cup_{\alpha' \in {\cal A} \setminus \{\alpha\}} V_{\alpha'})$.
As in (a) above it follows that $F_0 = V_{\alpha} \dot{+}   V'_{\alpha}$.
Therefore, this net is strongly convergent on the subspace $F_0 \dot{+}   {\hat R}$
to the projection $p_{\alpha}$ which has the properties (i) and (ii).

Finally, assume that $W_{\alpha}$ is a Coxeter group of type $B_{\Delta}$. Then
for any finite subset $\sigma \subset \Delta$ the product $u_{\sigma} =
\prod_{\delta \in \sigma} s_{\delta}$ of pairwise commuting reflexions
$s_{\delta} = s_{\epsilon_{\delta}, \epsilon_{\delta}^*} \in {\rm IR}_{\alpha}$
is an isometric involution with the fixed point subspace $\cap_{\delta \in
\sigma} {\rm Ker}\,\epsilon_{\delta}^* \supset {\hat R} \dot{+}   V'_{\alpha}$.
Similarly as above, the net of the restrictions $(u_{\sigma}\,|\,(F_0 \dot{+}
{\hat R}))$ is strongly convergent to an isometric involution $u_{\alpha}$
which has $ V_{\alpha}$ and ${\hat R} \dot{+}   V'_{\alpha}$ as its spectral
subspaces. It is easily seen that the projection $p_{\alpha} = (1_{F_0 \dot{+}
{\hat R}} + u_{\alpha})/2$ possesses all the properties mentioned in (i), (ii)
and (iii).

\noin $c$. By the closed graph theorem it is enough to check the second
statement. For a finite subset $\sigma \subset {\rm IR}_{\cal A}$ let
$W_{\sigma}$ be a finite group generated by reflexions from $\sigma$, and let
$V_{\sigma}$ be the linear span of the reflexion vectors of these reflexions.
Then the action of $W_{\sigma}$ in $V_{\sigma}$ is fixed point free, and so the
projection $p'_{\sigma} = (1/{\rm card}\,W_{\sigma})\sum _{g \in W_{\sigma}} g$
onto the fixed point subspace $F_{\sigma} \supset {\hat R}$ of $W_{\sigma}$
vanishes on $V_{\sigma}$. Consider the net of finite dimensional projections
$(p_{\sigma} = 1_{F_0 \dot{+}   {\hat R}} - p'_{\sigma}\,|\,(F_0 \dot{+}   {\hat R}))$
onto the subspaces $V_{\sigma}$. Observe that $\bigcup_{\sigma} V_{\sigma}$ is
dense in the subspace $F_0$. Since all of $p_{\sigma}$ vanish on ${\hat R}$ and
satisfy the norm inequalities of (ii), this net is strongly convergent to the
projection $p_{F_0 , {\hat R}}$, which also satisfies these inequalities. This
completes the proof. $\,\,\,\bigcirc$

\bigskip

\noin {\it Remark.} For further information on the Hilbert decomposition, see
Proposition 6.2 and examples 6.8 below.

\bigskip

\noin {\bf 5.7. Theorem.}
{\sl The subspaces $F, F_0, \hat R$ and $R_0$ are invariant with respect to the
group ${\rm Iso}\,E$, and there are natural monomorphisms ${\rm Iso}\,E
\hookrightarrow {\rm Iso}\,R_0 \times {\rm Iso}\,F_0$ , $G_0(E) \hookrightarrow
\prod\limits_{\beta \in \cal B} G_0(H_{\beta})$ and $\prod\limits_{\beta \in
\cal B} {\rm O} (H_{\beta}) \hookrightarrow {\rm Iso}\,(F \dot{+}   R_0)$.}

\bigskip

\noin {\sl Proof.} The invariance of the subspaces $F$ and $F_0$ was already
established in 5.2; the invariance of $R_0$ follows from the remark in 5.3.
Similar arguments applied to the conjugate action of ${\rm Iso}\,E$ on $E^*$
provide the invariance of $\hat R$.

Since the set $\{e \in S(E)\,|\,s_{e, e^*} \in IR(E)\}$ is contained in
$F_0\cup(\bigcup\limits_{\beta \in \cal B} H_{\beta}) \subset R \dot{+}   F_0$, the
latter summands being invariant, it follows from Lemma 4.1 that the restriction
mappings $g \longmapsto g\,|\,R_0\,,\, g \longmapsto g\,|\,F_0\,,\,g
\longmapsto g\,|\,H_{\beta}$ induce the monomorphisms ${\rm Iso}\,E
\hookrightarrow {\rm Iso}\,R_0 \times {\rm Iso}\,F_0$ and $G_0(E)
\hookrightarrow \prod\limits_{\beta \in \cal B} G_0(H_{\beta})$.

As for the last statement, fix arbitrary $g = \prod\limits_{\beta \in \cal
B}\bar u_{\beta} \in \prod\limits_{\beta \in \cal B} {\rm O} (H_{\beta}) $. For
any finite subset $\sigma \subset \cal B$ put $u_{\sigma} = \prod\limits_{\beta
\in \sigma} u_{\beta}$, where $u_{\beta} = \bar u_{\beta} p_{\beta} +(1_E -
p_{\beta}) \in {\rm Iso}\,E$ (see Theorem 1.b). We will show that the net
$\{u_{\sigma}\,|\,(F \dot{+}   R_0)\} \subset {\rm Iso}\,(F \dot{+}   R_0)$ strongly
converges to an element $u \in {\rm Iso}\,(F \dot{+}   R_0)$ such that
$u\,|\,H_{\beta} = \bar u_{\beta}$. Therefore, the correspondence
$\prod\limits_{\beta \in \cal B} {\rm O} (H_{\beta}) \owns g \longmapsto u \in
{\rm Iso}\,(F \dot{+}   R_0)$ yields the desired monomorphism.

By the Banach - Steinhaus Theorem, it is enough to show that for any $x \in F
\dot{+}   R_0$ the generalized sequence $(u_{\sigma}(x))$ is convergent. Let $x =
x_1 + x_2$, where $x_1 \in F$ and $x_2 \in R_0$. For any $\epsilon > 0$ there
exists a finite subset $\sigma \subset \cal B$ such that $||(1_E - \sum
\limits_{\beta \in  \sigma} p_{\beta})(x_2)||_E < \epsilon /2$. If $\sigma'$
and $\sigma''$ are two finite subsets of $\cal B$ containing $\sigma$, then
$u_{\sigma'} - u_{\sigma''} = (u_{\sigma'} - u_{\sigma''})(1_E -
\sum\limits_{\beta \in  \sigma} p_{\beta})$, and so $$||(u_{\sigma'} -
u_{\sigma''})(x_2)||_E \le ||u_{\sigma'}(1_E - \sum\limits_{\beta \in  \sigma}
p_{\beta})(x_2)||_E + ||u_{\sigma''}(1_E - \sum\limits_{\beta \in  \sigma}
p_{\beta})(x_2)||_E < \epsilon\,\,\,.$$ Thus, $(u_{\sigma}(x))$ is a
generalized Cauchy sequence, and hence it is convergent. This proves the
theorem. $\,\,\,\bigcirc$

\bigskip

\noin {\it Remark.} In general, the monomorphisms in Theorem 5.7 are not
surjective; see Example 6.8.2 below.

\section{ An application: Isometry groups of ideal generalized sequence spaces}

\noin {\bf 6.1.} {\it  Definitions.} Recall the following notions (see e.g.
[17], [19]). Let $(e_{\alpha})_{\alpha\in \Delta}$ be a system of vectors in a
Banach space $E_0$. It is called {\it a generalized Shauder basis} of $E_0$ if
each vector $e \in E_0$ has a unique, up to permutations, decomposition $e =
\sum_{i=1}^{\infty} a_i e_{\alpha_i}$, where $(\alpha_i)_{i=1,\dots}$ is a
sequence of pairwise distinct indicies from $\Delta$. If this series is still
convergent to $e$ after any permutation of its members, then this basis is
called {\it unconditional}. In this case for any choices of signes $\theta =
(\theta_{\alpha})_{\alpha \in \Delta}$, where $\theta_{\alpha}=\pm 1$, the
linear operators $M_{\theta}(e) = \sum_{i=1}^{\infty}\theta_{\alpha_i} a_i
e_{\alpha_i}$ are uniformly bounded. The number ${\rm sup}_{\theta}
\,||M_{\theta}||_{E_0}$ is called {\it the unconditional constant} of the basis
$(e_{\alpha})_{\alpha\in \Delta}$. For instance, any complete orthonormal
system in a Hilbert space is an unconditional basis with the unconditional
constant $1$. If the index set $\Delta$ is countable, we have the usual notion
of an unconditional basis.

The generalized unconditional basis $(e_{\alpha})_{\alpha\in \Delta}$ is called
{\it symmetric} if for any bijection $\pi: \Delta \rightarrow \Delta$ the
linear operator $$\pi^* : E_0 \owns e = \sum_{i=1}^{\infty} a_i e_{\alpha_i}
\mapsto  \sum_{i=1}^{\infty} a_i e_{\pi(\alpha_i)} = \pi^* (e) \in E_0$$ is
bounded, and so the infinite symmetric group $S_{\Delta} = {\rm
Biject}(\Delta)$ acts in $E_0$, being uniformly bounded there. The constant
${\rm sup}_{\theta, \pi} ||M_{\theta}\pi^*||_{E_0}$ is called {\it the
symmetric constant} of the basis
$(e_{\alpha})_{\alpha\in \Delta}$.

For instance, in the classical Banach space $c_0 (\Delta)$ of generalized
sequences convergent to zero, with $\Delta$ as a set of indicies, the system of
the standard basis vectors $(\epsilon_{\delta})_{\delta\in\Delta}$ form a
symmetric basis with the symmetric constant $1$ (remark that each vector in
$c_0 (\Delta)$ has a countable support). Fixing a generalized unconditional
basis $(e_{\alpha})_{\alpha\in \Delta}$ in $E_0$ we obtain a representation of
$E_0$ as a generalized sequence space contained in $c_0 (\Delta)$. If the
unconditional constant of this basis is $1$, then $E_0$ is an ideal Banach
lattice.

Recall that {\it an ideal generalized sequence space} $E$ is a Banach space of
sequences defined on an index set  $\Delta$ such that if
$x=(x_{\alpha})_{\alpha \in \Delta} \in E$, then for any sequence $y =
(y_{\alpha})_{\alpha \in \Delta}$ with $|y_{\alpha}| \le |x_{\alpha}|$ for all
$\alpha \in \Delta$ one has $y \in E$ and $||y||_E \le ||x||_E$. It is called
{\it a symmetric generalized sequence space} if $E$ is an ideal generalized
sequence space, where the symmetric group $S_{\Delta}$ of all bijections of
$\Delta$ acts isometrically. \\

The next simple lemma should be well known; by the lack of references we give a
proof. We say that a family of reflexions is {\it orthogonal} if the reflexions
from the family pairwise commute.

 \bigskip

\noin {\bf 6.2. Lemma.} {\sl Let $E$ be a Banach space with a total orthogonal
family of isometric reflexions
$(s_{\delta}=s_{\epsilon_{\delta},\epsilon_{\delta}^*})_{\delta \in \Delta}$.
Identify $E$ with a generalized sequence space with the index set $\Delta$ by
posing ${\bar x}= (\epsilon_{\delta}^*(x))_{\delta \in \Delta}$ for $x \in E$.
Let $E_0 = {\overline {\rm span}}\, (\epsilon_{\delta}\,|\,\delta \in \Delta)$.
Then we have

\smallskip

\noin a) The system $(\epsilon_{\delta})_{\delta \in \Delta}$ is a generalized
unconditional basis in $E_0$ with the unconditional constant $1$, and so $E_0$
is an ideal generalized sequence space.

\smallskip

\noin b) If the Coxeter group $B_{\Delta}$ of permutations and sign changes of
finite number of coordinates acts isometrically in $E_0$, then
$(\epsilon_{\delta})_{\delta \in \Delta}$ is a symmetric basis in $E_0$ with
the symmetric constant $1$, and so $E_0$ is a symmetric generalized sequence
space.}

\bigskip

\noin {\sl Proof. a}. Let $\sigma$ be a finite subset of $\Delta$. Consider the
coordinate subspace $$E_{\sigma} = \{x=(x_{\delta})_{\delta\in\Delta}\in
E_0\,|\,x_{\delta} =0 \,\,{\rm for\,\, all}\,\,\delta \notin \sigma\}\,\,\,.$$
Let $p_{\sigma} = {1\over 2} (1_{E_0} - u_{\sigma} )$, where $u_{\sigma} =
\prod_{\delta \in \sigma} s_{\delta}$, be the coordinate projection $E_0 \to
E_{\sigma}$. Since $u_{\sigma} \in {\rm Iso}\, E_0$, we have $||p_{\sigma}
||_{E_0} = ||1_{E_0} - p_{\sigma}||_{E_0} = 1$.

Fix an arbitrary vector $x \in E_0$. For any $n \in {\rm \bf N}$ there exists a
finite subset $\sigma_n$ of $\Delta$ and $y_n \in E_{\sigma_n}$ such that $||x
- y_n||_E < 1/n$. Then also $||p_{\sigma_n}(x) - y||_{E_0} < 1/n$, and so
$||(1_{E_0} - p_{\sigma_n})(x)||_{E_0} < 2/n$.  It follows that $x$ has at most
countable support contained in $\Omega = \bigcup_{i=1}^{\infty} \sigma_i =
\{\delta_1,\dots, \delta_k,\dots\}$, and $||x - \sum_{i=1}^k
x_i\epsilon_{\delta_i}||_{E_0} \to 0$. Thus the system
$(\epsilon_{\delta})_{\delta \in \Delta}$ is a generalized Shauder basis in
$E_0$. It is easily seen that for any fixed subset $\Omega \subset \Delta$ the
net of isometric involutions $(u_{\sigma}\,|\,\sigma \subset \Omega,\,{\rm
card}\,\sigma < \infty)$ strongly converges on $E_0$ to the isometric
involution $u_{\Omega}$, and therefore the basis $(\epsilon_{\delta})_{\delta
\in \Delta}$ of $E_0$ is unconditional with the unconditional constant $1$.

\smallskip

\noin $b$. Fix a permutation $\pi \in S_{\Delta}$, a vector $x \in E_0$ and
$\epsilon > 0$ arbitrarily. Let $\sigma$ be a finite subset of $\Delta$ such
that $||(1_{E_0} - p_{\sigma})(x)||_{E_0} < \epsilon$. There exists a finite
permutation $\pi' \in S_{\Delta}$ such that $\pi' | \sigma = \pi | \sigma$.
Since the Coxeter group $B_{\Delta}$ acts isometrically on $E_0$, we have
$$||\pi^* p_{\sigma}(x)||_{E_0} = ||\pi'^*p_{\sigma}(x)||_{E_0} =
||p_{\sigma}(x)||_{E_0}\,\,\,.$$ Thus, the linear operator $\pi^*$ is well
defined and isometric on the dense subspace ${\bf R}^{\Delta}$ of $E_0$.
Therefore, it can be extended isometrically onto $E_0$, and since the same is
true for $(\pi^{-1} )^*$, this extension does belong to the group ${\rm
Iso}\,E_0$. This proves  ($b$).$\,\,\,\bigcirc$

\bigskip

\noin {\it Remark.} It is not true, in general, that under the assumption of
this lemma $E$ itself should be an ideal space if all single sign changes are
isometries of $E$. As an example, consider the space $c$ of convergent
sequencies, which is not an ideal lattice.

\bigskip

We return to the Hilbert decomposition, keeping all the notation and the
conventions of Section 5.\\

\bigskip

\noin {\bf 6.3. Proposition.} {\it There exists an ideal generalized sequence
space $X$ with $\cal B$ as an index set such that the subspace $R_0$ of $E$ is
isometric to the Banach sum $(\bigoplus\limits_{\beta \in \cal B}
H_{\beta})_X$.}

\bigskip

\noin {\sl Proof.} For each $\beta \in \cal B$ fix a vector $e_{\beta} \in
S(H_{\beta})$. Consider the subspace $X = {\overline {\rm span}}\, (e_{\beta}
\,|\,\beta \in {\cal B}) \subset F \dot{+}   R_0$. Since the system of functionals
$(e_{\beta}^*\,|\,\beta \in {\cal B})$ is biorthogonal to the system
$(e_{\beta} \,|\,\beta \in {\cal B})$ and the reflexions $s_{e_{\beta},
e_{\beta}^*}\,|\,X$ are isometric, by Lemma 5.9.a, the latter system is an
unconditional basis in $X$ with the unconditional constant $1$, and so $X$ can
be identified with an ideal generalized sequence space on $\cal B$.

Put $R_1 = (\bigoplus\limits_{\beta \in \cal B} H_{\beta})_X$. We will show
that  the correspondence $$\tau : R_0 \owns x \longmapsto (p_{\beta}
(x))_{\beta \in {\cal B}} \in R_1$$ is a linear isometry of $R_0$ onto $R_1$.
Put $e'_{\beta} = {p_{\beta} (x) \over ||p_{\beta} (x)||_E} \in S(H_{\beta})$
if $p_{\beta} (x) \neq 0$. Let $\bar u_{\beta} \in {\rm O}(H_{\beta})$ is such
that $\bar u_{\beta}(e'_{\beta}) = e_{\beta}$ if $p_{\beta} (x) \neq 0$ and
$\bar u_{\beta} = 1_{H_{\beta}}$ otherwise. As follows from Theorem 5.7, there
exists $u \in {\rm Iso}\,R_0$ such that $u\,|\,H_{\beta} = \bar u_{\beta}$ for
all $\beta \in \cal B$. Since $u(x) \in X$ and $u$ is an isometry, it is clear
that $\tau (x) \in R_1$ and $||\tau (x)||_{R_1} = ||x||_{R_0}$.

To show that $\tau$ is surjective, fix arbitrary vector $\bar x \in R_1\,,\,
\bar x = (x_{\beta} \in H_{\beta})_{\beta \in \cal B}$. Then $x' =
\sum\limits_{\beta \in \cal B} ||x_{\beta}||_E\, e_{\beta} \in X \subset R_0$.
For each $\beta \in \cal B$ let $\bar u_{\beta} \in {\rm O}(H_{\beta})$ be such
that $\bar u_{\beta}(x_{\beta}) = ||x_{\beta}||_E\, e_{\beta}$. As above, there
exists $u \in {\rm Iso}\,R_0$ such that $u\,|\,H_{\beta} = \bar u_{\beta}$ for
all $\beta \in \cal B$. If $x_0 = u^{-1} (x')$, then we have $\tau (x_0) = \bar
x$. Thus, $\tau$ is an inversible isometry. This completes the proof.
$\,\,\,\bigcirc$

\bigskip

\noin {\bf 6.4.} {\it Notation.} Consider again an ideal generalized sequence
space $E$ with an index set $\Delta$. Without lost of generality one may assume
that $||\epsilon_{\delta}||_E = 1$ for all $\delta \in \Delta$. For a subset
$\Omega \subset \Delta$ let $E({\Omega})$ be a strip $E({\Omega}) =
\{x=(x_{\delta} ) \in E\,|\,x_{\delta} = 0$ for all $\delta \in \Delta
\setminus \Omega\}$. Any such strip is biorthogonally complemented in $E$;
indeed, the operator of multiplication by the characteristic function of
$\Omega$ is a bicontractive projection $p_{\Omega}\,:\,E \to E(\Omega)$ with
$1_E - p_{\Omega} = p_{\Delta \setminus \Omega}$.

{Put $E_0 (\Omega) = {\overline {\rm span}}\,(\epsilon_{\delta})_{\delta \in
\Omega}$, so that $E = E(\Delta), \,E_0 = E_0 (\Delta)$ and $E_0 (\Omega) = E_0
\cap E (\Omega)$. We also preserve in this particular case all the other
notation introduced in section 5. The next proposition shows that the Hilbert
and Coxeter decompositions of an ideal generalized sequence space yield an
orthogonal decomposition into strips.

A reflexion vector $\epsilon_{\delta}$ of the single sign change $s_{\delta} =
s_{\epsilon_{\delta}, \epsilon_{\delta}^*} \in {\rm IR} (E)$ belongs to a
certain subspace $V_{\alpha}$ or $H_{\beta}$. Putting $\Delta_{\alpha} =
\{\delta \in \Delta\,|\,s_{\delta} \in V_{\alpha}\}$ and $\Delta_{\beta} =
\{\delta \in \Delta\,|\,s_{\delta} \in H_{\beta}\}$ we obtain a disjoint
partition of $\Delta$ by the subsets $\{\Delta_{\alpha}
\,,\,\Delta_{\beta}\}_{\alpha \in {\cal A}\,,\,\beta \in {\cal B}}$. Put also
$\Delta_{{\cal A}} = \bigcup\limits_{\alpha \in {\cal A}} \Delta_{\alpha}$ and
$\Delta_{{\cal B}} = \bigcup\limits_{\beta \in {\cal B}} \Delta_{\beta}$, so
that $\Delta = \Delta_{{\cal A}} \cup \Delta_{{\cal B}}$.

\bigskip

\noin {\bf Proposition 6.5.} {\sl In the notation as above one has

\smallskip

\noin a)

\noin i) $V_{\alpha} = E_0 (\Delta_{\alpha})\,,\,{\hat V_{\alpha}} =
E(\Delta_{\alpha})$ for each $\alpha \in {\cal A}$ and

\noin ii) $H_{\beta} = E (\Delta_{\beta})=E_0 (\Delta_{\beta}) = l_2
(\Delta_{\beta})$ for each $\beta \in {\cal B}$;

\noin iii) if ${\rm card}\,(\Delta_{\alpha}) = \infty$, then $W_{\alpha}$ is a
Coxeter group of type $B_{\Delta_{\alpha}}$;

\smallskip

\noin b)

\noin i) $F_0 = E_0 (\Delta_{{\cal A}})$ and $R_0 = E_0 (\Delta_{{\cal B}})$,
so that $E_0 = R_0 \dot{+}   F_0$;

\noin ii) $F = E(\Delta_{{\cal A}})$ and ${\hat R} = E(\Delta_{{\cal B}})$, so
that $E = {\hat R} \dot{+}   F$;

\smallskip

\noin c)

\noin i) $p_{\alpha} = p_{\Delta_{\alpha}}\,|\,(F_0 \dot{+}   {\hat R})$ and
$p_{\beta} = p_{\Delta_{\beta}}$ for all $\alpha \in {\cal A}\,,\,\beta \in
{\cal B}$;

\noin ii) $p_{R_0, F} = p_{\Delta_{\cal B}}\,|\,(R_0 \dot{+}   F)$ and $p_{F_0,
{\hat R}} = p_{\Delta_{\cal A}}\,|\,(F_0 \dot{+}   {\hat R})$ (see Proposition
5.6).}

\bigskip

\noin {\sl Proof. a}. Put $\cal C = \cal A \cup \cal B$ and $V_{\gamma} =
H_{\gamma}$ for $\gamma \in \cal B$. By the above definitions, $E_0
(\Delta_{\gamma}) \subset V_{\gamma}$ for all $\gamma \in \cal C$. Let $\gamma
\in \cal C$, $\delta
\in \Delta_{\gamma}$ and $\delta' \in \Delta \setminus \Delta_{\gamma}$. By
Lemma 5.4, $s_{\delta'} = s_{\epsilon_{\delta'}, \epsilon_{\delta'}^*} \in {\rm
IR} (E)\setminus IR_{\gamma}$ commutes with any reflexion $s \in IR_{\gamma}$,
and so $V_{\gamma} \subset {\rm Ker}\,\epsilon_{\delta'}^*$. Therefore,
${\hat V_{\gamma}} \subset {\rm Ker}\,\epsilon_{\delta'}^*$, too, and hence
${\hat V_{\gamma}} \subset \bigcap\limits_{\delta' \in \Delta \setminus
\Delta_{\gamma}} {\rm Ker}\,\epsilon_{\delta'}^* = E(\Delta_{\gamma})$. In
particular, each reflexion vector $e$ of a reflexion $s=s_{e, e^*} \in IR(E)$
belongs to one of the strips $E(\Delta_{\gamma})$, where $\gamma \in \cal C$.
Namely, $e = (x_{\delta}) \in E(\Delta_{\gamma})$ iff $x_{\delta}\neq 0$ for at
least one $\delta \in \Delta_{\gamma}$. Furthermore, in the latter case either
$e = \pm \epsilon_{\delta}$ or the reflexions $s$ and $s_{\delta}$ do not
commute.

Let $\gamma = \alpha \in \cal A$. From the classification of the infinite
Coxeter groups in section 4 it follows that if $W$ is such a group and $s \in
W$, then the set of all reflexions in $W$ that do not commute with $s$ contains
not more than a finite subset of pairwise commuting reflexions. This means that
the reflexion vector $e$ of any given reflexion $s\in IR_{\alpha}$ has only a
finite number of non-zero coordinates, i.e. $e \in {\rm
span}\,(\epsilon_{\delta}\,|\,\delta \in \Delta_{\alpha}) \subset E_0
(\Delta_{\alpha})$. Thus, $V_{\alpha} \subset E_0 (\Delta_{\alpha})$, and
therefore,   $V_{\alpha} = E_0 (\Delta_{\alpha})$, which is the first statement
of ($a.i$). In particular, the reflexion vectors $(\epsilon_{\delta})_{\delta
\in \Delta_{\alpha}}$ of sign change reflexions $(s_{\delta})_{\delta \in
\Delta_{\alpha}} \subset IR_{\alpha}$ form a complete orthogonal system in
$V_{\alpha}$.

If ${\rm card}\,(\Delta_{\alpha}) <\infty$, then clearly $E
(\Delta_{\alpha})=E_0 (\Delta_{\alpha})=V_{\alpha}=\hat {V_{\alpha}}$. If ${\rm
card}\,(\Delta_{\alpha}) =\infty$, then by Corollary 4.6, the group
$W_{\alpha}$ generated by reflexions from $IR_{\alpha}$ is a Coxeter group of
type $A_{\Delta'}$ or $B_{\Delta'}$. But the Coxeter group $A_{\Delta'}$ does
not contain a complete set of pairwise commuting reflexions, i.e. there is no
orhtogonal subsystem of the root system $(\epsilon_{\delta} -
\epsilon_{\delta'}\,|\,\delta, \delta' \in \Delta'\,,\,\delta \neq \delta')$
which would be a Hamel basis of ${\bf R}_0^{\Delta'}$. This excludes the first
case, and so the group $W_{\alpha}$ should be a Coxeter group of type
$B_{\Delta'}$. It is clear that ${\rm card}\,(\Delta') = {\rm
card}\,(\Delta_{\alpha})$. This proves ($a.iii$).

Let, further, $\gamma = \beta \in \cal B$. Then $E_0 (\Delta_{\beta})$ is a
subspace of the Hilbert space $H_{\beta}$, and the system
$(\epsilon_{\delta})_{\delta \in \Delta_{\beta}}$ is an orthonormal basis of
$E_0 (\Delta_{\beta})$. Thus, $E_0 (\Delta_{\beta})= {\rm l}_2
(\Delta_{\beta})$. Assume that $H_{\beta} \neq E_0 (\Delta_{\beta})$ . Let $x
\in H_{\beta}$ be a non-zero vector orthogonal to $E_0 (\Delta_{\beta})$. It is
easily seen that $\epsilon_{\delta}^* (x) = 0$ for all $\delta \in
\Delta_{\beta}$. This is impossible, since $H_{\beta}\subset E
(\Delta_{\beta})$ and $x \neq \bar 0$. Therefore, $H_{\beta} = E_0
(\Delta_{\beta}) = {\rm l}_2 (\Delta_{\beta})$.

Let $p_{\beta} \,:\,E \to H_{\beta}$ be the projection as in Theorem 1.b.
Suppose that $H_{\beta}\neq E (\Delta_{\beta})$. Then the restriction
$p_{\beta}\,|\,E (\Delta_{\beta})$ is a non-identical projection, so that there
exists a non-zero vector $x \in {\rm Ker}\,p_{\beta}\cap E (\Delta_{\beta})$.
Fixing $\delta \in \Delta_{\beta}$, consider the plane $L = {\rm span} (x,
\epsilon_{\delta})$. There are two commuting isometric reflexions in $L$,
namely
$(1_L -2p_{\beta})\,|\,L$ and $s_{\delta}\,|\,L$. Therefore, $x \in {\rm
Ker}\,e_{\delta}^*$ for all $\delta \in \Delta_{\beta}$, and so $x = \bar 0$,
which is a contradiction. Hence, $H_{\beta} = E(\Delta_{\beta}) = E_0
(\Delta_{\beta}) = {\rm l}_2 (\Delta_{\beta})$. This proves ($a$), besides the
second equality in ($a.i$), which is proven below.

\smallskip

\noin $b$. For any $\gamma \in {\cal C}$ consider the isometric involution
$u_{\gamma} = 1_E - 2p_{\Delta_{\gamma}}$ with the spectral subspaces
$E(\Delta_{\gamma})$ and $E(\Delta \setminus \Delta_{\gamma})$. It is easily
seen that for any $s=s_{e, e^*} \in IR(E)$ the isometries $s u_{\gamma}$ and
$u_{\gamma} s$ coincide on the total system of reflexion vectors
$(\epsilon_{\delta})_{\delta \in \Delta}$. From Lemma 4.1 it follows that they
coincide on $E$.  Thus, the involution $u_{\gamma}$ commutes with each
reflexion $s_{e, e^*} \in IR(E)$. Therefore, one of its spectral subspaces
contains the vector $e$ and another one is contained in the mirror hyperplane
${\rm Ker}\, e^*$. Hence, for any $s_{e, e^*} \in IR_{\gamma}$ one has ${\rm
Ker}\, e^* \supset E(\Delta \setminus \Delta_{\gamma})$.

Let the set $\cal C$ be devided into two disjoint parts ${\cal C} = {\cal C}'
\cup {\cal C}''$. Put $\Omega' = \bigcup\limits_{\gamma \in {\cal C}'}
\Delta_{\gamma}\,,\,\Omega'' = \bigcup\limits_{\gamma \in {\cal C}''}
\Delta_{\gamma}$, so that $\Omega', \Omega''$ consist of some parts of the
disjoint partition $\Delta = \bigcup\limits_{\gamma \in {\cal C}}
\Delta_{\gamma}$. We are going to show, more generally, that ${\widehat {E_0
(\Omega')}} = E(\Omega')$, which easily implies the equalities in ($b.ii$) and
($a.i$).

By the considerations above, we have $E(\Omega') \subset \bigcap ({\rm
Ker}\,e^*\,|\,s_{e, e^*} \in {\rm IR}_{{\cal C}''})$, where ${\rm IR}_{{\cal
C}''} = \bigcup\limits_{\gamma \in {\cal C}''} {\rm IR}_{\gamma}$. On the other
hand, $E(\Omega') = \bigcap\limits_{\delta \in \Omega''} {\rm
Ker}\,\epsilon_{\delta}^* \supset \bigcap ({\rm Ker}\,e^*\,|\,s_{e, e^*} \in
IR_{{\cal C}''})$. Therefore, $E(\Omega') =  \bigcap ({\rm Ker}\,e^*\,|\,s_{e,
e^*} \in {\rm IR}_{{\cal C}''}) = {\widehat {E_0 (\Omega')}}$. The last
equality is clear from the definition of the envelope, because $s_{e, e^*} \in
IR_{{\cal C}''}$ iff $E_0 (\Omega') \subset {\rm Ker}\,e^*$. This proves ($b$)
and the second equality in ($a.i$).

\smallskip

$c$.  The isometric involutions $u_{\beta} = 1_E - 2p_{\Delta_{\beta}}$ and
$1_E - 2p_{\beta}$ coincide on vectors of the system
$(\epsilon_{\delta})_{\delta \in \Delta}$, so by Lemma 4.1, they coincide on
$E$.
This proves the second equality in ($c.i$). By the same reasoning (see
Proposition 5.6.$b.iii$) the first equality in ($c.i$) holds. The equalities
($c.ii$) follow from ($b$), just by the definition of the projections involved.
By Proposition 5.6.$b.i$, the projection $p_{\alpha}$ commutes with any sign
change reflexion $s_{\delta}$. By the same type of arguments as those used in
the proof of ($b$), it follows that ${\rm Ker}\, p_{\alpha} = {\rm Ker}
\,(p_{\Delta_{\alpha}}\,|\,(F_0 \dot{+}   {\hat R}))$. Since the images also
coincide,
we have the first equality in ($c.i$). This proves the proposition.
$\,\,\,\bigcirc$

\bigskip

This proposition, together with Theorem 5.7 and the remark that the union of
the subspaces $H_{\beta}\,,\,\beta \in \cal B,$ is invariant with respect to
the group ${\rm Iso}\,E$, leads to the following

\bigskip

\noin {\bf 6.6. Corollary.}
{\it
a) $$G_0(E) \subset \bigoplus\limits_{\beta \in {\cal B}} G_0 (l_2
(\Delta_{\beta}))\,\,\,\,\,\, and \,\,\,\,\,\,{\rm Iso}\,E \subset {\rm Iso}\,E
(\Delta_{{\cal A}}) \bigoplus {\rm Iso}\,E (\Delta_{{\cal B}})\,\,;$$
b) each element of the group $({\rm Iso}\,E)\,|\,E (\Delta_{\cal B})$ is of the
form $(x_{\beta}) \mapsto (u_{\beta} (x_{\beta}))$, where $x_{\beta} \in {\rm
l}_2 (\Delta_{\beta})$, $\,\,u_{\beta} :  {\rm l}_2 (\Delta_{\beta}) \to  {\rm
l}_2 (\Delta_{\pi (\beta)})$ is an isometry of Hilbert spaces for each $\beta
\in {\cal B}$, and $\pi$ is a permutation of the set $\cal B$.}

\bigskip

\noin {\bf 6.7.} {\it Remark.} Let $\alpha \in \cal A$ be such that ${\rm
card}\, \Delta_{\alpha} = \infty$. Then ${\rm Iso}\,V_{\alpha}$ contains a
Coxeter subgroup of type $B_{\Delta_{\alpha}}$. It is not true in general that
it contains also the symmetric group $S_{\Delta_{\alpha}}$ of shift operators.
In fact, this latter group is contained in ${\rm Iso}\,V_{\alpha}$, but
probably in some other representation of $V_{\alpha}$ as an ideal generalized
sequence space. Indeed, consider any symmetric generalized sequence space $M$
on $\Delta_{\alpha}$, such that the system  $(\epsilon_{\delta})_{\delta \in
\Delta}$ is a Shauder basis in $M$. Fix a disjoint partition of
$\Delta_{\alpha}$ into pairs $(\delta, \delta')$. Then by Lemma 6.2.$a$, the
corresponding subsystem $(\epsilon_{\delta} \pm \epsilon_{\delta'})$ of the
root system is an unconditional Shauder basis in $M$ with the unconditional
constant $1$, and this basis is not symmetric.  Thus, using the dual system of
functionals, one can represent the strip component $M = V_{\alpha}$ as an ideal
generalized sequence space which is not symmetric and such that the isometry
group does not act as permutations and sign changes (the image of a basis
vector under an isometric reflexion might be a vector with 4 non-zero
coordinates!). Recall that a symmetric basis is unique; moreover, a basis which
in a sense is symmetric enough, is unique [11].  Thus, here we have an
unconditional basis with a relatively small group of symmetries.

Conversely, if $g \in {\rm Iso}\, E$ is such that $g(V_{\alpha}) =
V_{\alpha'}$, where $\alpha \in \cal A$ is as above, then one can represent
$V_{\alpha}$ resp. $V_{\alpha'}$ as a symmetric generalized sequence space on
$\Delta_{\alpha}$ resp. $\Delta_{\alpha'}$, and then $g$ should be an operator
of the form $(x_{\alpha}) \mapsto  (\pm x_{\pi(\alpha)})$, where
$\pi\,:\,\Delta_{\alpha} \to \Delta_{\alpha'}$ is a bijection (indeed, $\pi$
must transfer the sign change reflexions from ${\rm IR}_{\alpha}$ into sign
change reflexions from ${\rm IR}_{\alpha'}$).

In [24], [25] certain conditions on an ideal generalized sequence space are
given which guarantee that its isometry group acts by permutations and sign
changes. This is always the case in a symmetric sequence spaces different from
${\rm l}_2$ [23, ch. IX], [6] (see also [2], [8] for the complex field).
\bigskip

Next we give several examples related to the results of Sect. 5 and 6.

\bigskip

\noin {\bf 6.8.} {\it Examples.}

\medskip

\noin 1) ([24], [25]) Fix a sequence of real numbers $p_k \ge 1\,,\,k=1,\dots$.
{\it The Orlicz--Nakano space} $E = l(\{p_k\})$ consists of all sequences of
real numbers $x = (\xi_k)_{k=1}^{\infty}$ such that the following norm is
finite: $$||x||_E = {\rm inf}\,\{\lambda > 0\,|\,\sum\limits_{k=1}^{\infty}
|\xi_k / \lambda |^{p_k} \le 1\}\,\,.$$
It is an ideal sequence space.

Put $\Delta_q = \{i \in {\bf N}\,|\,p_i = q\}$, where $q= q_1 , q_2 ,\dots$ are
pairwise distinct. Then ${\cal A} = \{q_i\,|\,i \neq 2\}$ and $\Delta_{\cal A}
= \{i\,|\,p_i \neq 2\}\,,\, {\cal B} = \{2\}$ if $\Delta_2 \neq \emptyset$ and
${\cal B} = \emptyset$ otherwise; $E(\Delta_q) = {\rm l}_q (\Delta_q)$. The
group ${\rm Iso}\, E$ is the direct product of the groups ${\rm O}\, ({\rm l}_2
(\Delta_2))$ and ${\rm Iso}\,\Delta_{\cal A}$, where ${\rm Iso}\,\Delta_{\cal
A}$ is the group of all permutations of coordinates $\xi_i,\,i \in \Delta_{\cal
A}$, preserving the partition $\Delta_{\cal A} = \bigcup \Delta_{q_i}$, and
arbitrary sign changes of these coordinates. Indeed, this direct product
evidently is a subgroup of ${\rm Iso}\, E$; the converse inclusion follows from
the results of [24], [25] in view of the decomposition from Corollary 6.6.

In a similar way one can describe the isometry groups of more general modular
sequence spaces or of Banach sums of (symmetric) ideal sequence spaces.

\medskip

\noin 2) Let $E$ be the space of all convergent complex sequences with the
supremum norm. Then $E$ is a Banach sum of the real euclidean planes
$H_i\,,\,i=1,\dots$. We have that $E = \hat R$, and $E_0 = R_0$ is the subspace
of sequences in $E$ convergent to zero. The group ${\rm Iso}\, E_0$ is the
semi-direct product of ${\rm O}(2)^{\omega}$ and the infinite symmetric group
$S_{\omega}$, while ${\rm Iso}\, E$ is its proper subgroup (indeed, if $g \in
{\rm Iso}\, E$, then the corresponding sequence of orthogonal plane
transformations from ${\rm O}(2)^{\omega}$ is convergent). This shows that all
the inclusions in Theorem 5.7 are strict. Observe that here $E$ is not an ideal
sequence space.

\medskip

\noin 3) Consider $E = {\bf R}^4$ with the norm
$$||x||_E = ||(\xi_1 , \xi_2, \varsigma_1, \varsigma_2)||_E =
[((\xi_1^2 +     \xi_2^2)^{1/2} + |\varsigma_1 |)^2 +
\varsigma_2^2]^{1/2}\,.$$
It is easily seen that here
$\hat R = R_0 = H = \{x \in E\,|\,\varsigma_1 = \varsigma_2 = 0\}$
and $F=F_0=\{x \in E\,|\,\xi_1 = \xi_2 =0\}$.
Furthermore, $E = F_0\dot{+}   R_0$ is an ideal space, and both strips $F_0$ and
$R_0$ are euclidean planes. Thus, $G_0 (E) \neq G_0 (F_0) \oplus G_0 (R_0)$,
and so ${\rm Iso}\, E \neq {\rm Iso}\, F_0 \oplus {\rm Iso}\, R_0$ (cf.
Corollary 6.6.$a$).

\medskip

\noin 4) Slightly modifying example 4, consider $E \cong {\bf R}^6$ with the
norm
$$||x||_E = ||(\xi_1 , \xi_2, \eta_1, \eta_2, \varsigma_1, \varsigma_2)||_E =
[((\xi_1^2 + \xi_2^2)^{1/2} + |\varsigma_1 |)^2 + ((\eta_1^2 + \eta_2^2)^{1/2}
+ |\varsigma_2|)^2]^{1/2}\,.$$
Being the direct sum of two  euclidean planes  $ H_1$ and $ H_2$, which are
strips invariant under $G_0 (E)$, the subspace $R_0$ itself is euclidean. Thus,
$G_0 (R_0) \neq G_0 (H_1) \oplus G_0 (H_2) = G_0 (E)$ (cf. Corollary 6.6.$a$).

\medskip

\noin 5) Let, further, $\bar E = {\rm c} \oplus {\rm c} \oplus {\rm c}$ with
the norm $||(x, y, z)||_{\bar E} = \sup_{i=1,\dots} \,\{(\xi_i^2 +
\eta_i^2)^{1/2} + |\varsigma_i |\}$, where $x = (\xi_i)_{i=1}^{\infty} \in {\rm
c}\,,\, y = (\eta_i)_{i=1}^{\infty} \in {\rm c}\,,\,z =
(\varsigma_i)_{i=1}^{\infty} \in {\rm c}$. Consider the hyperplane $E = \{(x,
y, z) \in \bar E\,|\,\lim_{i \to \infty} \eta_i = \lim_{i \to \infty}
\varsigma_i\}$. Here we have ${\hat R} \approx {\rm c} \oplus {\rm c}_0\,,\,F
\approx {\rm c}_0$.  Thus, $ {\hat R} \dot{+}   F$ is a hyperplane in $E$, and
there is no contractive projection of $E$ onto ${\hat R}$ and onto $F$,  in
contrary to the case of ideal sequence spaces (cf. Propositions 5.6 and
6.5.$b,\,c$).

\bigskip

The following questions are directly related to the subject of this paper.

\smallskip

For a given Banach space $E$, consider the constant $$c(E) = \inf_{e\in E,
\,e^* \in E^*,\,e^* (e) = 1}\{||s_{e, e^*}||_E\}\,\,.$$ It is clear that $1 \le
c(E) \le 3$, and $c(E) = 1$ in the case when there exists an isometric
reflexion is $E$. It is easily seen that $c(L_p )$ is a convex function of $p$
which takes the value $1$ only for $p=2$ and the value $3$ only for $p=1$ and
$p=\infty$. For any given finite set of reflexions in E one can find an
equivalent norm $||\cdot||'$ on $E$ in such a way that the group generated by
these reflexions will be a subgroup of the isometry group of the new norm. In
particular, $c(E, ||\cdot||') = 1$.

Consider, further, the constant $$\varsigma(E) = \sup_{||\cdot||' \sim
||\cdot||_E} \{c(E, ||\cdot||')\}\,\,.$$ By the definition, $\varsigma (E) \in
[1,\,3]$ is a numerical invariant of isomorphism. Is it nontrivial?

Let $\varsigma (n) = \varsigma ({\bf R}^n) = \varsigma ({\rm l}^n_2)$. Let
$M_n$ be the Minkowski compact of classes of isometric norms in ${\bf R}^n$,
endowed with the Banach--Mazur distance. Denote by $A_n$ the subset of $M_n$
which consists of the classes of norms having an isometric reflexion (or, the
same, a hyperplane of symmetry). It is easy to show that $\log\varsigma (n)$
coincides with the radius of the metric factor space $M_n / A_n$ with respect
to the distinguish point which corresponds to $A_n$. It is known that the
radius of the $M_2$ centred at the class of euclidean norms is $\log\sqrt{2}$
(F. Behrend, 1937; see [12, sect. 7] for this and for some further
information). Thus, $\varsigma (2) \le \sqrt{2}$.

\medskip

\noin {\bf 6.9.} {\it Problem.} Is it true that\\
$\varsigma (3) < 3$ ? \\
$\varsigma (n) < 3 $ for any $n$ ? \\
$\limsup_{n \to \infty} \varsigma (n) < 3$ ?\\
$\varsigma ({\rm l}_2) < 3$ ?

\medskip

\noin If the answer to any of the above questions is ``yes'', which seems to be
less plausible, then, of course, the exact value of the corresponding constant
$\varsigma$ would be worthwhile to find.

\bigskip


\begin{center} {\LARGE References} \end{center}

\noin [1] Yu. Abramovich, M. Zaidenberg. {\em  A rearrangement invariant space
isometric to $L_p$ coincides with $L_p$},
Interaction between Functional Analysis, Harmonic Analysis, and
Probability, N. Kalton, E. Saab eds., Lect. Notes in Pure Appl. Math. 175,
Marcel Dekker, New York, 1995, 13-18.

\noin [2] J. Arazy, {\em Isometries of complex symmetric sequence spaces},
Math. Z., 188 (1985), 427--431.

\noin [3] S. Banach, {\em Th\'eorie des op\'erations lin\'eaires}, Hafner
Publishing Co, N.Y. 1932.

\noin [4] J.J. Barry, {\em On the convergence of ordered sets of projections},
Proc. Amer. Math. Soc., 5 (1954), 313--314.

\noin [5] N. Bourbaki, {\em Groupes et alg\`ebres de Lie}, Ch.IV - VI, Hermann,
Paris, 1968.

\noin [6] M. Sh. Braverman, E.M. Semenov, {\em Isometries of symmetric spaces},
          Soviet Math. Doklady, 15 (1974), 1027--1030.

\noin [7] M.M. Day, {\em Normed linear spaces}, Springer, N.Y. e.a., 1973.

\noin [8] R.J. Fleming, J.E. Jamison, {\em Isometries on certain
Banach spaces}, J. London Math. Soc., 2 (1974),  121--127.

\noin [9] R. Godement, {\em Th\'eor\'emes taub\'eriens et th\'eorie spectrale},
Ann. Sci. Ecole Norm. Sup., (3) 64 (1947), 119--138.

\noin [10] Y. Gordon, D. R. Lewis, {\em Isometries of diagonally symmetric
Banach spaces}, Israel J. of Mathem., 28 (1977), 45--67.

\noin [11] Y. Gordon, R. Loewy, {\em Uniqueness of $(\Delta)$ bases and
isometries of Banach spaces}, Math. Ann., 241 (1979), 159--180.

\noin [12] B. Gr\"unbaum, {\em Measures of symmetry for convex sets}, Proc.
Symp. Pure. Mathem. VII Convexity (1963), AMS, Providence, R.I., 233--270.

\noin [13] N.J. Kalton, B. Randrianantoanina, {\em Isometries on
rearrangement-invariant spaces}, C. R. Acad. Sci. Paris, 316 (1993), 351--355.

\noin [14] N.J. Kalton, B. Randrianantoanina, {\em  Surjective isometries on
rearrangement-invariant spaces}, Quart. J. Math. Oxford, (2) 45 (1994),
301--327.

\noin [15] N. J. Kalton, G.W. Wood, {\em  Orthonormal systems in Banach spaces
and their applications}, Math. Proc. Cambr. Phil. Soc., 79 (1976), 493--510.

\noin [16] M. A. Krasnosel'skii, {\em On a spectral property of compact linear
operators in the space of continuous functions}, Problems of the mathematical
analysis of complex systems, 2, Voronezh (1968),  68--71 (Russian).

\noin [17] S. G. Krein, Ju. I. Petunin, E. M. Semenov. {\em Interpolation of
linear operators}, Providence. R.I., 1982.

\noin [18] P.--K. Lin, {\em Elementary isometries of rearrangement--invariant
spaces}, preprint, 1995, 1--22.

\noin [19] J. Lindenstrauss and L. Tzafriri, {\em Classical Banach Spaces,
Vol.1,2}, Springer, 1977, 1979.

\noin [20] Yu.I. Lyubich, {\em  On the boundary spectrum of a contraction in
Minkovsky spaces}, Sib. Mat.Zh., 11:2 (1970), 358--369.

\noin [21] Yu.I. Lyubich, L. N. Vaserstein, {\em Isometric embeddings between
classical Banach spaces, curbature formulas, and spherical designs}, Geometriae
Dedicata, 47 (1993), 327--362.

\noin [22] B. Randrianantoanina, {\em Isometric classification of norms in
rearrangement--invariant function spaces}, preprint, 1995, 1--15.

\noin [23] S. Rolewicz, {\em  Metric Linear Spaces}, Polish Scientific
Publishers, Warsaw, 1972.

\noin [24] A.I. Skorik, {\em Isometries of ideal coordinate spaces}, Uspekhi
Mathem. Nauk, 31:2 (1976), 229--230 (Russian).

\noin [25] A.I. Skorik, {\em On isometries of a class of ideal coordinate
spaces}, Teorija Funkziy, Funkz. Analyz Prilozh., 34 (1980), 120--131
(Russian).

\noin [26] A.I. Skorik, M.G. Zaidenberg, {\em Groups of isometries containing
reflexions}, Function. Anal. Appl., 10 (1976), 322--323.

\noin [27] A.I. Skorik, M.G. Zaidenberg, {\em Groups of isometries containing
reflexions}, preprint VINITI, DEP N\, 1638-78 (1978), 43p. (Russian; English
transl.: {\em On isometric reflexions in Banach spaces},
Pr\'epublication de l'Institut Fourier des Math\'ematiques, 267, Grenoble 1994,
36p.

\noin [28] J. Wermer, {\em  The existence of invariant subspaces}, Duke
Math.J.,           19 (1952), 615--622.

\noin [29] M.G. Zaidenberg, {\em Groups of isometries of Orlich spaces},
Soviet  Math. Dokl., 17 (1976), No.2, 432--436.

\noin [30] M.G. Zaidenberg, {\em On the isometric classification of symmetric
spaces}, Soviet Math. Dokl., 18 (1977), 636--640.

\noin [31]  M.G. Zaidenberg, {\em Special representations of isometries of
functional spaces}, Investigations on the theory of functions of several real
variables, Yaroslavl' 1980 (Russian; English transl.: {\em A representation of
isometries on function spaces}, Pr\'epublication
de l'Institut Fourier des Math\'ematiques, 305, Grenoble 1995, 7p.).

\bigskip

\noin M. Zaidenberg, A. Skorik

\noin Institut Fourier des Math\'ematiques

\noin Universit\'e Grenoble I

\noin BP 74

\noin 38402 Saint Martin d'H\`eres-c\'edex

\noin France

\medskip

\noin E-mail: zaidenbe@puccini.ujf-grenoble.fr

\end{document}